\documentclass{nmeauth}

\usepackage{cite}
\usepackage{amsmath}
\usepackage{amsfonts}
\usepackage{tikz}
\usepackage{subfigure}
\usepackage{textcomp}

\newcommand{\vect}[1]{\underline{#1}}
\newcommand{\matr}[1]{\underline{\underline{#1}}}
\newcommand{\transpose}[1]{{#1}^{\text{T}}}

\newcommand{\SIJ}{S_{IJ}}
\newcommand{\DeltauIJ}{\vect{\Delta u}_{IJ}}
\newcommand{\nIJ}{\vect{n}_{IJ}}
\newcommand{\sIJ}{\vect{s}_{IJ}}
\newcommand{\tIJ}{\vect{t}_{IJ}}
\newcommand{\DIJ}{D_{IJ}}
\newcommand{\DIJeq}{D_{IJ}^{0}}
\newcommand{\PIJ}{P_{IJ}}

\newcommand{\FIJ}{\vect{F}_{IJ}}
\newcommand{\MIJ}{\vect{M}_{IJ}}

\newcommand{\epsilonIJ}{\varepsilon_{IJ}^v}

\newcommand{\gradvect}{\vect{\nabla}}
\newcommand{\gradmatr}{\matr{\nabla}}
\newcommand{\divscal}{\text{div}\text{~}}
\newcommand{\divvect}{\vect{\text{div}}\text{~}}
\newcommand{\rot}{\vect{\text{curl}}\text{~}}
\newcommand{\trace}{\text{tr~}}

\DeclareTextSymbol{\degre}{T1}{6}
\DeclareTextSymbol{\degre}{0T1}{23}

\newtheorem{thm}{Theorem}[section]
\newtheorem{lem}[thm]{Lemma}

\RequirePackage[dvips]{graphicx}
\RequirePackage{times}

\begin{document}

\NME{0}{0}{0}{0}{10}

\runningheads{L.~Monasse, C.~Mariotti}{An energy-preserving D.E.M. for elastodynamics}
\noreceived{}
\norevised{}
\noaccepted{}
\title{An energy-preserving Discrete Element Method for elastodynamics}

\author{Laurent Monasse\affil{1}\comma\corrauth and Christian Mariotti\affil{2}}

\address{\affilnum{1}\ Universit\'e Paris-Est, CERMICS, 6 et 8 avenue Blaise
  Pascal, Cit\'e Descartes -- Champs-sur-Marne, 77455 Marne-la-Vall\'ee
  Cedex 2, France\\
\affilnum{2}\ CEA DAM DIF, F-91297 Arpajon, France}

\corraddr{Laurent Monasse, Universit\'e Paris-Est, CERMICS, 6 et 8 avenue Blaise
  Pascal, Cit\'e Descartes -- Champs-sur-Marne, 77455 Marne-la-Vall\'ee
  Cedex 2, France}

\begin{abstract}
  We develop a Discrete Element Method (DEM) for elastodynamics using
  polyhedral 
  elements. We show that for a given choice of forces and torques, we
  recover the equations of linear elastodynamics in small
  deformations. Furthermore, the torques and
  forces derive from a potential energy, and thus the global equation is
  an Hamiltonian dynamics. The use of an explicit symplectic time
  integration scheme allows us to recover conservation of energy, and
  thus stability over long time simulations. These theoretical results
  are illustrated by numerical simulations of test cases involving large
  displacements. 
\end{abstract}


\section{Introduction}

Particle methods are meshless simulation techniques in which a continuum medium
is approximated through the dynamics of a set of interacting
particles. Two main classes of particle methods can be distinguished :
Discrete Element methods (DEM), which rely on the contact interaction of
material particles by means of forces and torques, and Smooth Particle
Hydrodynamics (SPH) methods, in which the continuum is discretized by
localized kernel functions. 

Discrete Element methods consist in the resolution of the equations of
motion of a set of particles submitted to forces and torques. It is thus
possible to account for a variety of phenomena (behaviour laws, models,
scales,...) using a single numerical method. A wide variety of Discrete
Elment methods have been designed changing the expression of the forces,
with particular attention devoted to specific aspects. Discrete Element methods
have first been developed by Hoover, Arhurst and Olness \cite{Hoover} in
models for crystalline materials. Their application to geotechnical problems was
carried out by Cundall and Strack \cite{Cundall}, and their use in
granular materials and rock simulation is still widespread
\cite{Wolf,ref5}. Discrete Element Methods have also been used to
simulate thermal conduction in granular assemblies \cite{ref4} or
fluid-structure interaction \cite{ref9}. The model is also able to
account for grain size effects \cite{ref61}, and to treat fracture in a natural
way. Discrete Element methods used for granular materials generally
describe particles as
spherical elements interacting via noncohesive, frictional contact forces
\cite{Wolf}.  For brittle materials, models also use unilateral contact
forces, combined with bonds which simulate cohesion \cite{ref5}. Kun and
Herrmann developed a combination of the contact model with a lattice
model of beams to account for the cohesion \cite{ref57}, which has been
extended to Reissner models of beams to simulate large rotations of the
material \cite{ref56,ref61}. The authors use Voronoi tesselations to
generate the polygonal particles. However, the results obtained still
depend on the size of the discretization (which physically corresponds to
the size of heterogeneities) \cite{ref61}. The effective macroscopic
Young modulus and Poisson ratio highly depend on the isotropy of the
distribution of the particles and are only empirically linked to their
microscopic value for the Reissner beams \cite{ref57}. 

In a different approach, SPH methods describe the particles as smooth
density kernel functions. The kernel functions are an approximation of
the partition of unity. The continuous equations of evolution of the
fluid or solid material therefore induce the dynamics of the particles.
Originating from astrophysical compressible fluid simulations
\cite{Monaghan1,Lucy}, SPH was extended to incompressible fluids
\cite{Monaghan2} and to elastic and plastic dynamics \cite{Libersky},
and used for fluid-structure interaction with both domains discretized
with SPH \cite{ref6}. A state of the art review of the method with
applications to solid mechanics is presented in
\cite{Hoover2}. SPH preserves the total mass of the
system exactly. However, in tensile regime, unphysical clusters of
particles tend to appear in situations where a homogeneous response is
expected \cite{ref58}. Hicks, Swegle and Attaway advocate the
smoothing of the variables between neighbouring particles to stabilize
the method, rather than introducing artificial viscosities
\cite{ref59}. Bonet and Lok have addressed the issue of angular
momentum preservation, and show that rotational invariance is
equivalent to the exact evaluation of the gradients of linear velocity
fields, which can be achieved either through correction of the kernel
function or through a modification of its gradient \cite{ref55}. 
In order to circumvent the difficulties affecting SPH, Yserentant
developed the Finite Mass method, in which particles of fixed size and
shape also possess a rotational degree of freedom (spin). The method achieves
effective partition of unity, and thus preserves momentum, angular
momentum and energy, ensuring stability \cite{ref60}.

The Moving Particle Semi-implicit (MPS) method is a variant of the SPH
method developed by Koshizuka. It consists in the derivation of the
dynamics of a set of points from a 
discrete Hamiltonian \cite{Koshizuka1}. As in the SPH method, the
differential operators are approximated by a kernel function of compact
support. The expression of the approximated differential operators is
inserted in the classical Hamiltonian of the system, and by application
of Hamilton's equations, the dynamics of the discretized system is
obtained. To preserve the Hamiltonian structure of the dynamic of the
system through time discretization, the authors use symplectic schemes
\cite{ref54}. The MPS method has been used initially for free-surface
flows \cite{Koshizuka1,Koshizuka2}, and has been extended to nonlinear
elastodynamics \cite{Koshizuka3,ref54} and to fluid-structure
interaction \cite{ref12}. Using similar ideas, by deriving the dynamics
of the system from a discrete Hamiltonian, Fahrenthold has simulated
compressible flows \cite{ref15} and impact events with breaking of the
target \cite{ref1,ref2}.

These methods show the importance of the preservation of
momentum and energy for the accuracy and stability of the scheme over
long-time simulation. The use of symplectic schemes ensures the
preservation of the structure of Hamilton's equations by the numerical
time integration, and therefore the preservation of momentum and energy
\cite{Hairer}. Simo, Tarnow and Wong note, however, that while ensuring
the stability of the simulation for small time steps, the symplectic
schemes fail to preserve exactly energy and become unstable for larger
time steps \cite{ref18}. They derive a general class of implicit
time-stepping algorithms which exactly enforce the conservation of
momentum, angular momentum and energy. The algorithms are built in order
to preserve linear and angular momentum, and energy conservation is
enforced either with a projection method (projection on the manifold of
constant energy) or with a collocation method. The algorithm is used for
nonlinear elasticity in large deformation using finite element methods
\cite{ref18,ref17,ref19} and for low-velocity impact
\cite{ref16}. 

In this article, we extend and analyze the Discrete Element method
initially introduced by Mariotti
\cite{Mariotti}. Combining a Discrete Element Method with a
lattice model of beams, we are able to account for the cohesion of the
material, and analytically recover the macroscopic behaviour of the
continuous material. The method, Mka3D, has been
successfully used to simulate the 
propagation of seismic waves in linear elastic medium \cite{Mariotti}.
Here, we extend the properties of the algorithm to the case of large
displacements without
fracture. Contrary to usual Discrete Element methods, we are able to
derive the microscale forces and torques analytically from the
macroscopic Young modulus and Poisson ratio, and to prove the
convergence of the method as the grid is refined. In addition, as in MPS
methods, we derive the forces and torques between particles from a
Hamiltonian formulation. Using a symplectic scheme, we ensure the
preservation of energy over long-time simulations, and thus stability of
the method. This allows for the simulation of three-dimensional wave
propagation as well as shell or multibody dynamics.
The paper is organized as follows.
In section \ref{sec:description}, we describe the lattice model used. We
introduce the Hamiltonian of the system and we derive the
expression of forces and torques chosen to simulate linear
elasticity. In section \ref{sec:precision}, we show that these
expressions lead to a macroscopic behaviour of the material equivalent
to a Cosserat continuum, with a characteristic length of the order of
the size of the particles. Hence, the model is consistent with
a Cauchy continuum medium up to second-order accuracy, in the case of small
displacement and small deformation. The microscopic values of Young
modulus and Poisson ratio yield directly the macroscopic values, and we
can choose Poisson ratio in the whole interval $(-1,0.5)$.
In section \ref{sec:time-scheme}, we then describe the
symplectic RATTLE time-scheme \cite{Hairer}, which allows us to
preserve a discrete energy over long-time simulations. These theoretical
results are illustrated by numerical simulations of test cases involving
large displacements in section \ref{sec:resultsnum}.

\section{Description of the method}

\label{sec:description}

\subsection{Geometrical description of the system}

In order to discretize the continuum material, several methods have
been suggested for Discrete Element Methods. Most authors working
on granular materials use hard 
spheres, in order to simplify the computation of contacts between
particles, as the exact form of the particles is mainly
unknown. However, in the case of the simulation of a continuous
material, this method is not adapted as the interstitial vacuum between
spheres is inconsistent with the compactness of the solid. In addition, the
difficulty to obtain a dense packing of hard spheres, and the problem of
the expression of cohesion between the particles, have led us to use
Voronoi tesselations instead, as suggested in \cite{ref57,ref56}. The
particles are therefore convex polyhedra which define a partition of the
entire domain. As we shall see, this method allows us to handle any
Poisson ratio $\nu$ strictly between $-1$ and $0.5$, independently from
the size of the particles. On the contrary, most granular sphere packing
methodologies account for a limited range of $\nu$, which is size dependent.

\begin{figure}
\centering
\resizebox{0.7\textwidth}{!}{
  \includegraphics{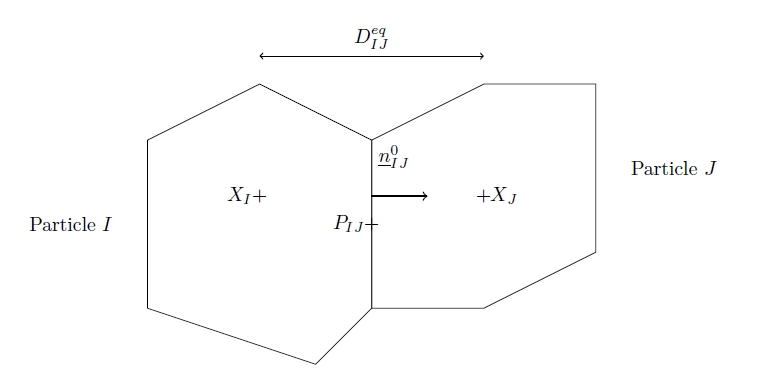}
}
\caption{Geometric description of the particles}\label{fig:particles}
\end{figure}

The following parameters are relevant to describe the motion
of a given particle $I$ : $\vect{X}_I$ and $\vect{v}_I$ denote
respectively the position and velocity of its center of mass
($\vect{v}_I = \frac{d\vect{X}_I}{dt}$), $\matr{Q}_I$ denotes the
orthogonal rotation matrix of the frame attached to the rigid particle,
and the angular velocity vector $\vect{\Omega}_I$ is uniquely defined by~:
  \begin{equation}
    \matr{j}(\vect{\Omega}_I) = \frac{d\matr{Q}_I}{dt}\matr{Q}_I^T, \label{eqn:deriveerotation}
  \end{equation}
  where the map $\matr{j} : \mathbb{R}^3
\rightarrow \mathbb{R}^{3\times3}$ is such that~:
\begin{equation*}
  \forall \vect{x}\in\mathbb{R}^3,\; \forall
  \vect{y}\in\mathbb{R}^3,\;
  \matr{j}(\vect{x})\cdot\vect{y}=\vect{x}\wedge\vect{y}
\end{equation*}
Finally, the material of particle $I$ is described by its mass $m_I$,
its volume $V_I$ and its principal moments of inertia $I_I^1$, $I_I^2$ 
and $I_I^3$. We suppose the local frame attached to the
particle is attached to the principal axes of inertia
$(\vect{e}^1_I,\vect{e}^2_I,\vect{e}^3_I)$.
The matrix of inertia in the fixed frame is given by~:
\begin{equation}
  \matr{R}_I = \matr{Q}_I \cdot \matr{R}_I^0 \cdot \matr{Q}_I^{-1} \label{eqn:inertie}
\end{equation}
where $\matr{R}_I^0$ is the matrix of inertia $\matr{R}_I^0$ written in
the inertial frame :
\begin{equation*}
  \matr{R}_I^0 = \left(
    \begin{array}{ccc}
      I_I^1 & 0 & 0 \\
      0 & I_I^2 & 0 \\
      0 & 0 & I_I^3
    \end{array} \right)
\end{equation*}
We also define the parameters $d_I^1$, $d_I^2$ and $d_I^3$ as :
\begin{equation*}
  d_I^i = \frac{I_I^1+I_I^2+I_I^3}{2} - I_I^i, \quad i=1,2,3
\end{equation*}
and we introduce the following matrix $\matr{D}_I$ defined in the
inertial frame~: 
\begin{equation*}
  \matr{D}_I = \left(
    \begin{array}{ccc}
      d_I^1 & 0 & 0 \\
      0 & d_I^2 & 0 \\
      0 & 0 & d_I^3
    \end{array} \right)
\end{equation*}

The Discrete Element Method relies on the computation of forces and
torques between nearest neighbours particles. We denote by $\mathcal{V}_I$
the list of the neighbouring particles linked to particle $I$. For each
link between two particles $I$ and $J$, we define $\PIJ$ the center of
mass of the interface, $\SIJ$ the surface of the interface, the
distance between particles $I$ and $J$~: 
\begin{equation*} 
  \DIJ = \lVert \vect{X_I X_J} \rVert, 
\end{equation*}
and the initial exterior normal vector for link $IJ$~:
\begin{equation*} 
  \nIJ=\frac{1}{\DIJ} \vect{X_I X_J}
\end{equation*}
We define two normalized orthogonal vectors of the interface $\sIJ$ and
$\tIJ=\nIJ \wedge\sIJ$, serving as
references to evaluate the torsion between particles $I$ and $J$.

These parameters are given a fixed value at the beginning of
the computation. $\DIJ^0$ and $\nIJ^0$ respectively denote the initial
values for $\DIJ$ and $\nIJ$. The particles are therefore assumed to be
rigid. However, compressibility effects are taken into account through
the expression of interaction potentials.

In addition, we define the following quantities~:
\begin{itemize}
\item the displacement at the interface between particles $I$ and $J$~:
  \begin{equation*}
    \DeltauIJ =  \vect{X_J} - \vect{X_I} +
    \matr{Q}_J\cdot\vect{X_J^0\PIJ} -
    \matr{Q}_I\cdot\vect{X_I^0\PIJ}
  \end{equation*}
\item When particle $I$ has several free interfaces (i.e. not linked
  to another particle), these surfaces are marked as stress-free. To
  account for the free deformation of the particle in these
  directions, free-volume $V_I^l$ is defined as the sum of the volumes
  of all pyramidal polyhedra with a free surface as basis and $X_I^0$
  as summit.  
\item the volumetric deformation $\varepsilon_I^v$ of particle $I$~ is
  defined as the sum of all contributions of the deformations of the
  material links of particle $I$. We have assumed that the
  bending of the link between two particles does not affect volume, as
  long as the centers of the interface of the two particles stay in
  contact. The corrective term on the volume is active only on particles
  having a free surface, and accounts for the boundary condition
  $\matr{\sigma}\cdot\vect{n}=\vect{0}$. We derive it in Appendix
  \ref{sec:freesurf}. 
  \begin{equation*}
    \varepsilon_I^v =
    \sum_{J\in\mathcal{V}_I}{\frac{1}{2}\frac{\SIJ}{V_I+3\frac{\nu}{1-2\nu}V_I^l}\DeltauIJ\cdot\nIJ}
  \end{equation*}
  \item The interpolated volumetric deformation for link $(IJ)$~:
    \begin{equation*}
      \epsilonIJ = \frac{1}{2}(\varepsilon_I^v+\varepsilon_J^v)
    \end{equation*}
\end{itemize}

\subsection{Expression of the Hamiltonian of the system}

We denote by $E$ the Young's modulus and by $\nu$ the Poisson's ratio for the
material.
The Hamiltonian formulation of the elastodynamic equations on a domain $\Omega$ is as follows~:
\begin{equation}
  H(\vect{q},\vect{p}) =
  \int_{\Omega}{\frac{1}{2\rho}\vect{p}\cdot\vect{p}} + U(\vect{q})
\end{equation}
where $\vect{q}$ is the displacement field and $\vect{p} = \rho\vect{v}$ is the
density of momentum. $U(\vect{q})$ is the potential energy
of the system. It can be expressed in terms of the stress tensor
$\matr{\sigma}$ and the linearized strain tensor $\matr{\varepsilon}=\frac{1}{2}(\matr{\nabla}\vect{q}+\transpose{\matr{\nabla}\vect{q}})$~:
\begin{equation}
  U(\vect{q}) = W(\matr{\varepsilon}) = \frac{1}{2}\int_{\Omega}{\matr{\sigma}(\matr{\varepsilon}):\matr{\varepsilon}}
\end{equation} 
In the case of Cauchy linear elasticity, we use the constitutive relation
\begin{equation}
  \matr{\sigma}(\matr{\varepsilon}) = \frac{E}{1+\nu}\matr{\varepsilon} + \frac{E\nu}{(1+\nu)(1-2\nu)}\trace(\matr{\varepsilon})\matr{Id}
\end{equation}
to derive the expressions of $W(\varepsilon)$ and $U(\vect{q})$~:
\begin{equation}
  W(\matr{\varepsilon}) =
  \frac{1}{2}\int_{\Omega}{\frac{E}{1+\nu}\matr{\varepsilon}:\matr{\varepsilon}+\frac{E\nu}{(1+\nu)(1-2\nu)}\trace(\matr{\varepsilon})^2} \label{eqn:ham_elastodynamique}
\end{equation}
\begin{equation}
  U(\vect{q}) =
  \frac{1}{2}\int_{\Omega}{\frac{E}{2(1+\nu)}\gradmatr\vect{q}:\gradmatr\vect{q}
  + \frac{E}{2(1+\nu)(1-2\nu)}(\divscal\vect{q})^2}
\end{equation}

We choose to discretize the Hamiltonian formulation as a discrete
Hamiltonian $H_h$. The displacement field $\vect{q}$ is derived from the
values of $(\vect{X}_I,\matr{Q}_I)$. The density of momentum derives from~:
\begin{align}
  \vect{T}_I &= m_I\vect{v}_I \\
  \matr{P}_I &= \matr{j}(\vect{\Omega}_I)\cdot\matr{Q}_I\cdot\matr{D}_I
\end{align}
We define~:
\begin{equation}
  H_h(\vect{X},\matr{Q},\vect{T},\matr{P}) =
  \frac{1}{2}\sum_I{\frac{1}{m_I}\vect{T}_I\cdot\vect{T}_I} + \frac{1}{2}\sum_I{\trace(\matr{P}_I\cdot\matr{D}_I^{-1}\cdot\transpose{\matr{P}_I})}
  + U_h(\vect{X},\matr{Q}) \label{eqn:hamiltonien}
\end{equation}

The discretized potential energy is split into three terms~:
\begin{equation*}
  U_h(\vect{X},\matr{Q}) = U_t(\vect{X},\matr{Q}) + U_d(\vect{X},\matr{Q}) + U_f(\matr{Q})
\end{equation*}
$U_t(\vect{X},\matr{Q})$ corresponds to the first term of
(\ref{eqn:ham_elastodynamique}) : we approach
the strain of the link $(IJ)$ in the direction $\nIJ$
$\matr{\varepsilon}\cdot\nIJ$ by the normalized displacement
$\frac{1}{\DIJeq}\DeltauIJ$, and we use the approximation~:
\begin{equation}
  \matr{\varepsilon}:\matr{\varepsilon} \approx \sum_{J\in\mathcal{V}_I}{(\matr{\varepsilon}\cdot\nIJ)^2}
\end{equation} 
We therefore write~:
\begin{equation*}
  U_t(\vect{X},\matr{Q}) = \frac{1}{2}\sum_{(IJ)}{\SIJ\frac{E}{1+\nu}\frac{\DeltauIJ\cdot\DeltauIJ}{\DIJeq}}
\end{equation*}
This energy accounts for the deformation of each link between two
particles.

$U_d(\vect{X},\matr{Q})$ corresponds to the second term of
(\ref{eqn:ham_elastodynamique}) : we approach the trace of the strain
$\text{tr}(\matr{\varepsilon})$ in particle $I$ by the sum of the
normalized displacements $\varepsilon_I^v$ for links surrounding $I$. A
corrective term is added for cells having a free boundary~:
\begin{equation*}
  U_d(\vect{X},\matr{Q})
  = \frac{1}{2}\sum_I{\frac{E\nu}{(1+\nu)(1-2\nu)}(V_I+3\frac{\nu}{1-2\nu}V_I^l)(\varepsilon_I^v)^2}
\end{equation*}
This energy accounts for the global volumetric deformation of each particle.

The former two terms are sufficient to recover the equations of
elastodynamics inside the solid. However, for the method to be able to
cope with thin one-element shells, we add the pure
flexion term $U_f(\matr{Q})$~: 
\begin{equation*}
    U_f(\matr{Q})
    = -\sum_{(IJ)}{\frac{\SIJ}{\DIJeq}\left(\alpha_n(\matr{Q}_{J}\cdot\nIJ^0)\cdot(\matr{Q}_{I}\cdot\nIJ^0)
    \right. }
    \left. + \alpha_s(\matr{Q}_{J}\cdot\sIJ)\cdot(\matr{Q}_{I}\cdot\sIJ) +
      \alpha_t(\matr{Q}_{J}\cdot\tIJ)\cdot(\matr{Q}_{I}\cdot\tIJ)
    \right)
  \end{equation*}
This term accounts for the flexion between particles. The coefficients
$\alpha_n$, $\alpha_s$ and $\alpha_t$ are chosen to recover the exact
flexion and torsion of a beam, and are detailed in Appendix \ref{sec:flexion}.

\subsection{Derivation of the forces and torques between particles}
\label{sec:derivation}

We use Hamilton's equations for the system (\ref{eqn:hamiltonien})~:
\begin{align}
  \vect{\dot{X}}_I &= \frac{\partial H_h}{\partial\vect{T}_I} \label{eqn:deplpos} \\
  \matr{\dot{Q}}_I &= \frac{\partial H_h}{\partial\matr{P}_I} \label{eqn:rotpos}\\
  \vect{\dot{T}}_I &= -\frac{\partial H_h}{\partial\vect{X}_I} \label{eqn:deplvit}\\
  \matr{\dot{P}}_I &= -\frac{\partial H_h}{\partial\matr{Q}_I} +
  \matr{\Lambda}_I\cdot\matr{Q}_I \label{eqn:rotvit}
\end{align}
where $\matr{\Lambda}_I$ is the symmetric matrix of the Lagrange
multipliers associated with the constraint
$\transpose{\matr{Q}_I}\cdot\matr{Q}_I = \matr{Id}$.

Equations (\ref{eqn:deplpos}) and (\ref{eqn:rotpos}) give us the usual
kinematic relations between position and velocity~:
\begin{align*}
  \vect{\dot{X}}_I &= m_I^{-1}\vect{T}_I = \vect{v}_I \\
  \matr{\dot{Q}}_I &= \matr{P}_I\cdot\matr{D}_I^{-1}  = \matr{j}(\vect{\Omega}_I)\cdot\matr{Q}_I 
\end{align*}

The derivation of forces and torques from the potential energies is
carried out in Appendix \ref{sec:forces}. We obtain
$m_I\vect{\dot{v}}_I=\FIJ$ where $\FIJ$, the force exerted by particle
$I$ on particle $J$, is given by~:
\begin{equation}
  \FIJ = \frac{\SIJ}{\DIJeq}\frac{E}{1+\nu}\DeltauIJ 
  + \SIJ\frac{E\nu}{(1+\nu)(1-2\nu)}\epsilonIJ\left(\nIJ+\frac{1}{\DIJ}\DeltauIJ-\frac{1}{\DIJ}(\DeltauIJ\cdot\nIJ)\nIJ\right)\label{eqn:forcetotale}
\end{equation}
This expression can be seen as a discrete version of Hooke's law of
linear elasticity
\begin{equation}
  \matr{\sigma} = \frac{E}{1+\nu}\matr{\varepsilon} +
  \frac{E\nu}{(1+\nu)(1-2\nu)}\text{tr}(\matr{\varepsilon})\matr{Id}
  \label{eqn:Hooke}
\end{equation}
using the previous analogies between $\frac{1}{\DIJeq}\DeltauIJ$ and
$\matr{\varepsilon}$, $\epsilonIJ$ and $\trace\matr{\varepsilon}$, and
noting that $\matr{\sigma}\cdot\vect{n}$ is a force per surface unit (a
pressure).

For the rotational part, we define the two following torques~:
\begin{equation}
  \MIJ^t = \frac{\SIJ}{\DIJeq}\frac{E}{1+\nu}
  (\matr{Q}_I\cdot\vect{X_I^0\PIJ})\wedge\DeltauIJ
  + \frac{E\nu}{(1+\nu)(1-2\nu)}\epsilonIJ\SIJ(\matr{Q}_I\cdot\vect{X_I^0\PIJ})\wedge\nIJ
  \label{eqn:torque}
\end{equation}
\begin{equation}
  \MIJ^f = \frac{\SIJ}{\DIJeq}\left(\alpha_n
  (\matr{Q}_I\cdot\nIJ^0)\wedge(\matr{Q}_J\cdot\nIJ^0)\right. 
  \left. + \alpha_s (\matr{Q}_I\cdot\sIJ)\wedge(\matr{Q}_J\cdot\sIJ)
  \right. 
  \left. + \alpha_t (\matr{Q}_I\cdot\tIJ)\wedge(\matr{Q}_J\cdot\tIJ)
  \right) \label{eqn:bendingtorsion}
\end{equation}
We note the fact that $\MIJ^t$ corresponds to the torque at the center
of mass of the force $\FIJ$
exerted by particle $J$ on particle $I$ at point $\PIJ$~:
\begin{equation*} 
  \MIJ^t = (\matr{Q}_I\cdot\vect{X_I^0\PIJ})\wedge\FIJ
\end{equation*}
and $\MIJ^f$ is the flexion-torsion torque.
We get the equation on the
angular velocity~:
\begin{equation}
  \frac{d}{dt}\left(\matr{R}_I\cdot\vect{\Omega}_I\right) =
\sum_{J\in\mathcal{V}_I}{\MIJ^t+\MIJ^f}\label{eqn:rfd2}
\end{equation}

In the case when exterior forces and torques are applied to the system,
they are to be added to the internal forces and torques computed
above.

\section{Consistency and accuracy of the scheme}
\label{sec:precision}

In this section, we investigate the consistency and the accuracy of the
scheme. We first propose a modified equation for small displacements and small
deformations. As the equations obtained are coupled dynamics for
displacement and rotation, we compare the model with Cosserat
generalized continuum, and recover a Cauchy continuum as the spatial
discretization $h$ tends to zero.

\subsection{Modified equation for the scheme}
\label{sec:equivalent}

The modified equation approach is a standard scheme analysis where a set
of continuous equations verified by the approximate solution is seeked
for. These modified equations should be an approximate version of
continuous equations derived from physics.

In order to be able to carry out a Taylor developments of the
displacement, we place 
the points of the Voronoi tesselation on a Cartesian grid. The Discrete
Element method can be seen, in this simplified case, as a Finite Difference
scheme. 

We assume that no exterior force and no exterior torque are applied on
the system. The displacement $\vect{\xi}_I$ of particle $I$ is given by~:
\begin{equation*}
  \vect{\xi}_I = \vect{X}_I - \vect{X}_I^0
\end{equation*}
We assume that $\vect{\xi}$ is a regular function on the domain, and
we can therefore expand $\vect{\xi}_J$ at point $I$ with Taylor
series if $J\in\mathcal{V}_I$. We denote $\Delta x$, $\Delta y$ and
$\Delta z$ the grid steps in each direction, and $h$ their maximum.

We assume displacements and rotations to be small. We denote
$\theta_x^I$, $\theta_y^I$ and $\theta_z^I$ the small rotation angles
around axes $x$, $y$ and $z$. 

Using (\ref{eqn:forcetotale}), a simple Taylor development of
the equations of motion yields for the displacement~:
\begin{multline}
  \rho \ddot{\xi}_x =
  \frac{E}{1+\nu}\left(\frac{\partial^2\xi_x}{\partial x^2} 
  + \frac{\partial^2\xi_x}{\partial y^2}
  + \frac{\partial^2\xi_x}{\partial z^2} 
  + \frac{\partial\theta_z}{\partial y} 
  - \frac{\partial\theta_y}{\partial z}\right) 
  + \frac{E\nu}{(1+\nu)(1-2\nu)}\left(\frac{\partial^2\xi_x}{\partial
  x^2}
  + \frac{\partial^2\xi_y}{\partial x\partial y} 
  + \frac{\partial^2\xi_z}{\partial x\partial z}\right) \\
  + \frac{E}{1+\nu}\left(\frac{\Delta x^2}{12}\frac{\partial^4
  \xi_x}{\partial x^4} 
  + \frac{\Delta y^2}{12}\frac{\partial^4 \xi_x}{\partial y^4} \right. 
  \left. + \frac{\Delta z^2}{12}\frac{\partial^4 \xi_x}{\partial z^4}  
  + \frac{\Delta y^2}{6}\frac{\partial^3\theta_z}{\partial y^3} 
  - \frac{\Delta z^2}{6}\frac{\partial^3\theta_y}{\partial z^3}\right) \\
  + \frac{E\nu}{(1+\nu)(1-2\nu)}\left(\frac{\Delta
  x^2}{3}\frac{\partial^4\xi_x}{\partial x^4} 
  + \frac{\Delta x^2}{6}\frac{\partial^4\xi_y}{\partial x^3\partial y} 
  \right.
  \left. + \frac{\Delta x^2}{6}\frac{\partial^4\xi_z}{\partial
  x^3\partial z} 
  + \frac{\Delta y^2}{6}\frac{\partial^4\xi_y}{\partial x\partial y^3} 
  + \frac{\Delta z^2}{6}\frac{\partial^4\xi_z}{\partial x\partial z^3}
  \right) \\
  + \mathcal{O}(h^3)  \label{eqn:devdeplacementx}
\end{multline}
The same results hold for $\xi_y$ and $\xi_z$ permuting the indices
$x$, $y$ and $z$ circularly.

Using (\ref{eqn:torque}) and
(\ref{eqn:bendingtorsion}), (\ref{eqn:rfd2}) gives the equivalent
equation for the rotation~:
\begin{multline}
  \frac{\Delta y^2+\Delta z^2}{12}\rho\ddot{\theta}_x =
  \frac{E}{1+\nu}\left(\frac{\partial\xi_z}{\partial y} 
  - \frac{\partial\xi_y}{\partial z} - 2\theta_x \right. 
  \left. + \frac{\Delta y^2}{6}\frac{\partial^3\xi_z}{\partial y^3} 
  - \frac{\Delta z^2}{6}\frac{\partial^3\xi_y}{\partial z^3} 
  + \frac{\Delta y^4}{120}\frac{\partial^5\xi_z}{\partial y^5} \right. \\
  \left. - \frac{\Delta z^4}{120}\frac{\partial^5\xi_y}{\partial z^5}  
  - \frac{\Delta y^2}{4}\frac{\partial^2\theta_x}{\partial y^2} 
  - \frac{\Delta z^2}{4}\frac{\partial^2\theta_x}{\partial z^2}\right. 
  \left. - \frac{\Delta y^4}{48}\frac{\partial^4\theta_x}{\partial y^4} 
  - \frac{\Delta z^4}{48}\frac{\partial^4\theta_x}{\partial z^4}\right)  \\
  + E\left[\frac{\Delta y^2+\Delta
  z^2}{12(1+\nu)}\left(\frac{\partial^2\theta_x}{\partial x^2} 
    + \frac{\Delta x^2}{12}\frac{\partial^4\theta_x}{\partial
  x^4}\right) \right. 
    \left. + \frac{\Delta z^2}{12}\left(\frac{\partial^2\theta_x}{\partial y^2}
    + \frac{\Delta y^2}{12}\frac{\partial^4\theta_x}{\partial y^4}\right)
    \right.  \\ 
    \left. + \frac{\Delta y^2}{12}\left(\frac{\partial^2\theta_x}{\partial z^2} 
    + \frac{\Delta z^2}{12}\frac{\partial^4\theta_x}{\partial z^4}\right)
    \right]  
  + \mathcal{O}(h^5) \label{eqn:devrotationx}
\end{multline}
The same results hold for $\theta_y$ and $\theta_z$ permuting the indices
$x$, $y$ and $z$ circularly.

We see that these sets of equations couple $\vect{\xi}$ and
$\vect{\theta}$, and by construction of the method, no constitutive
law exists between $\vect{\xi}$ and $\vect{\theta}$. The fact that a
rotation remains in the equations can be compared to Cosserat continuum
theory. We investigate this comparison in the following subsection.

\subsection{Comparison with Cosserat and Cauchy continuum theories}
\label{sec:cosserat}

In a Cosserat model for continuum media, the kinematics is described by
a displacement field $\vect{u}$ and a rotation field $\vect{\phi}$. A
modified strain tensor $\matr{\varepsilon}$ and a new curvature strain
tensor $\matr{\kappa}$ are introduced \cite{Eringen}~: 
\begin{align*}
  \matr{\varepsilon} &= \gradmatr\vect{u} + \matr{j}(\vect{\phi}) \\
  \matr{\kappa} &= \gradmatr\vect{\phi}
\end{align*}
We define $\matr{t}$ and $\matr{\mu}$ the stress and couple stress
tensors. We assume the following constitutive relations~:
\begin{align}
  \matr{t} &= \lambda \trace(\matr{\varepsilon})\matr{Id} +
  \mu\matr{\varepsilon} + \mu_c \matr{\varepsilon}^{\text{T}} \label{eqn:consteq1} \\
  \matr{\mu} &= \alpha\trace(\matr{\kappa})\matr{Id} + \gamma\matr{\kappa} +
  \beta\matr{\kappa}^{\text{T}} \label{eqn:consteq2}
\end{align}
where $\lambda$, $\mu$, $\mu_c$, $\alpha$, $\beta$ and $\gamma$ are
elastic moduli.

The dynamical equations for the system are~:
\begin{align*} 
\rho \vect{\ddot{u}} &= \divvect\matr{t}\\
\matr{I}_c \vect{\ddot{\phi}} &= \divvect\matr{\mu} + \vect{\matr{e}}:\matr{t}
\end{align*}
where $\rho$ denotes the density, $\matr{I}_c$ is a characteristic
inertia matrix, $:$ denotes the double contraction product of tensors,
and $\vect{\matr{e}}$ is defined as follows~:
\begin{equation*}
(\vect{\matr{e}})_{ijk} = \left\{
\begin{array}{ll}
1 & \text{if $(ijk)$ is an even permutation}\\
-1 & \text{if $(ijk)$ is an odd permutation}\\
0 & \text{otherwise}
\end{array}
\right.
\end{equation*} 

Using the constitutive relations (\ref{eqn:consteq1}) and
(\ref{eqn:consteq2}), the following equations can be obtained~:
\begin{align}
  \rho \vect{\ddot{u}} &= (\lambda+\mu_c)\gradvect\divscal\vect{u} +
  \mu\Delta\vect{u} + (\mu-\mu_c)\rot\vect{\phi}
  \label{eqn:deplcosserat} \\
  \matr{I}_c \vect{\ddot{\phi}} &=
  (\alpha+\beta)\gradvect\divscal\vect{\phi} + \gamma\Delta\vect{\phi}
  - 2(\mu-\mu_c)\vect{\phi} + (\mu-\mu_c)\rot\vect{u} \label{eqn:rotcosserat}
\end{align}

Identifying the terms of (\ref{eqn:deplcosserat}) with equation (\ref{eqn:devdeplacementx}), we find~:
\begin{align*}
  \lambda &= \frac{E\nu}{(1+\nu)(1-2\nu)} \\
  \mu &= \frac{E}{1+\nu}\\
  \mu_c &= 0
\end{align*}
and we therefore recover the classical expression, for Cauchy media, of
the first Lam\'e coefficient $\lambda_{Cauchy}$, and 
$\frac{\mu+\mu_c}{2}$ corresponds to the classical second Lam\'e coefficient
$\mu_{Cauchy}$.
Comparing then equation (\ref{eqn:rotcosserat}) with equation
(\ref{eqn:devrotationx}), we find~: 
\begin{equation*}
  \matr{I}_c = \rho\left(
    \begin{array}{ccc}
      \frac{\Delta y^2+\Delta z^2}{12} & 0 & 0 \\
      0 & \frac{\Delta x^2+\Delta z^2}{12} & 0 \\
      0 & 0 & \frac{\Delta x^2+\Delta y^2}{12}
    \end{array}
  \right)
\end{equation*}
For a given $h=\Delta x=\Delta y=\Delta z$, we see that the modified equations for
the scheme are those of a Cosserat generalized continuum, with
second-order accuracy, and the coefficients verify $\alpha+\beta = 0$
and $\gamma = \frac{E}{2(1+\nu)}h^2$. In the case of an anisotropic
mesh size ($\Delta x \neq\Delta y\neq \Delta z$), we cannot identify
the coefficients with the isotropic Cosserat equations, due to the
presence of the Laplacian operator. We can however find an anisotropic
Cosserat model with weighted second derivatives instead of the Laplacian.

One of the main characteristics of a Cosserat generalized continuum is
to exhibit a characteristic length for the material, $l_c$, which
describes the length of the nonlocal interactions. $l_c$ is defined as~:
\begin{equation*}
  l_c^2 = \frac{\gamma}{\mu+\mu_c}
\end{equation*}
In our case, we see that~:
\begin{equation*}
  l_c = \frac{\sqrt{2}}{2}h
\end{equation*}
$l_c$ is of the same
order as the size of the particles. In an
homogenization analysis framework, S. Forest,
F. Pradel and K. Sab have shown \cite{ref49} that when the macroscopic
length of the system is fixed and the characteristic length $l_c$ of
the Cosserat continuum tends to 0, the macroscopic behavior of the
material is that of a Cauchy continuum. We therefore converge to a
Cauchy continuum as $h$ tends to 0.

As a consequence, displacement $\vect{\xi}$, acceleration
$\vect{\ddot{\xi}}$, rotation $\vect{\theta}$ and acceleration
of rotation $\vect{\ddot{\theta}}$ in equations
(\ref{eqn:devdeplacementx}) and (\ref{eqn:devrotationx}) converge to
finite macroscopic quantities. Therefore, using the equations on
rotation, we find~:
\begin{equation}
  \vect{\theta} = \frac{1}{2}\rot\vect{\xi} + \mathcal{O}(h^2) \label{eqn:rotxi}
\end{equation}
which is the classical definition of the local rotation of a Cauchy
material at order 2. Using this relation in the equations of
displacement, we find the equations of linear elasticity for a Cauchy
continuum medium up to error terms of order $\mathcal{O}(h^2)$~:
\begin{equation*}
  \rho\vect{\ddot{\xi}} = \frac{E}{2(1+\nu)}\Delta\vect{\xi} +
  \frac{E\nu}{(1+\nu)(1-2\nu)}\gradvect\divscal\vect{\xi} 
  + \mathcal{O}(h^2)
\end{equation*}
and taking $\frac{1}{2}\rot$ of this equation, we find the equivalent
equation on rotation up to error terms of order $\mathcal{O}(h^2)$~:
\begin{equation}
  \rho\vect{\ddot{\theta}} = \frac{E}{2(1+\nu)}\Delta\vect{\theta} +
  \mathcal{O}(h^2) 
\end{equation}
We recover a second-order accuracy on the rotation $\vect{\theta}$. As
equation (\ref{eqn:rotxi}) shows, $\vect{\theta}$ is a derivate of
$\vect{\xi}$, and we should expect only first-order accuracy using a
second-order accurate method on $\vect{\xi}$. We have
therefore improved the accuracy on $\vect{\theta}$ using the Discrete
Element method.

\section{Preservation of the Hamiltonian structure by the time integration scheme}
\label{sec:time-scheme}

\subsection{Description of the scheme}

\label{sec:scheme}

The model built has a Hamiltonian 
structure. To preserve this property after time discretization,
we use a symplectic time integration scheme. As the system
(\ref{eqn:deplpos})--(\ref{eqn:rotvit}) is a constrained Hamiltonian
system \cite[Sec VII.5]{Hairer}, it is natural to use the following
RATTLE scheme \cite{Andersen} with time-step $\Delta t$~:  

\begin{align}
  \vect{T}_I^{n+1/2} &= \vect{T}_I^{n} - \frac{\Delta t}{2}\frac{\partial U_h}{\partial\vect{X}_I}(\vect{X}^{n}, \matr{Q}^{n}) \label{eqn:translation1} \\
  \matr{P}_I^{n+1/2} &= \matr{P}_I^{n} - \frac{\Delta
    t}{2}\frac{\partial U_h}{\partial\matr{Q}_I}(\vect{X}^{n}, \matr{Q}^{n}) 
  + \frac{\Delta t}{2}\matr{\Lambda}_I^n\matr{Q}_I^{n}
  \label{eqn:rotation1}
\end{align}
\begin{align}
  \vect{X}_I^{n+1} &= \vect{X}_I^{n} + \frac{\Delta
  t}{m_I}\vect{T}_I^{n+1/2} \label{eqn:translation2} \\
  \matr{Q}_I^{n+1} &= \matr{Q}_I^{n} +
  \Delta t\matr{P}_I^{n+1/2}\matr{D}_I^{-1}\label{eqn:rotation2}
\end{align}
\begin{equation}
  \text{where } \matr{\Lambda}_I^n \text{ is such that }\transpose{\matr{Q}_I^{n+1}}\cdot\matr{Q}_I^{n+1} =
  \matr{Id} \label{eqn:constraint1}  
\end{equation}
\begin{align}
  \vect{T}_I^{n+1} &= \vect{T}_I^{n+1/2} - \frac{\Delta
    t}{2}\frac{\partial U_h}{\partial\vect{X}_I}(\vect{X}^{n+1}, \matr{Q}^{n+1}) \label{eqn:translation3}\\
  \matr{P}_I^{n+1} &= \matr{P}_I^{n+1/2} - \frac{\Delta
    t}{2}\frac{\partial U_h}{\partial \matr{Q}_I}(\vect{X}^{n+1}, \matr{Q}^{n+1}) 
  + \frac{\Delta t}{2}\matr{\tilde{\Lambda}}_I^{n+1}\matr{Q}_I^{n+1}
  \label{eqn:rotation3} , 
\end{align}
\begin{equation}
  \text{where } \matr{\tilde{\Lambda}}_I^{n+1} \text{ is such that }\transpose{\matr{Q}_I^{n+1}}\cdot\matr{P}_I^{n+1}\cdot\matr{D}_I^{-1} +
\matr{D}_I^{-1}\cdot\transpose{\matr{P}_I^{n+1}}\cdot\matr{Q}_I^{n+1} =
\matr{0} \label{eqn:constraint2}
\end{equation}
where $\matr{\Lambda}_I^n$ and $\matr{\tilde{\Lambda}}_I^n$ are symmetric
matrices, the Lagrange multipliers associated with the constraints
(\ref{eqn:constraint1}) and (\ref{eqn:constraint2}).
We denote the scheme (\ref{eqn:translation1})--(\ref{eqn:constraint2})
by~:
\begin{equation*}
  (\vect{X}^{n+1},\matr{Q}^{n+1},\vect{T}^{n+1},\matr{P}^{n+1}) = \Psi_{\Delta t}(\vect{X}^{n},\matr{Q}^{n},\vect{T}^{n},\matr{P}^{n}) 
\end{equation*}

The proof for RATTLE's symplecticity can be found in
\cite{Leimkuhler}. As a consequence, in the absence of
exterior forces, the energy of the system is an invariant of the system,
and is preserved by the numerical integration in time. More precisely,
the error is of order $\mathcal{O}(e^{-\frac{\kappa}{\Delta t}})$ 
over a time period of $e^{\frac{\kappa}{\Delta t}}$, with $\kappa>0$
independent from $\Delta t$ \cite{Hairer}. This yields
the stability of the simulation over long time periods if the time step
is chosen sufficiently small. In addition, we directly derive from
(\ref{eqn:translation1})--(\ref{eqn:constraint2}) that the linear and
angular momentum are exactly preserved.

Another important property of the RATTLE scheme is its
reversibility. Starting with the knowledge of positions and velocities at
time $(n+1)\Delta t$, we recover the positions and velocities at time
$n\Delta t$ with the following scheme~:
\begin{equation*}
  (\vect{Q}_{T,n},\matr{Q}_{R,n},\vect{P}_{T,n},\matr{P}_{R,n}) = \Psi_{-\Delta t}(\vect{Q}_{T,n+1},\matr{Q}_{R,n+1},\vect{P}_{T,n+1},\matr{P}_{R,n+1}) 
\end{equation*}
As a reversible scheme, RATTLE is of even order, and as it is
consistent, it is a second-order scheme. 

RATTLE has the advantage of enforcing explicitly matrix $\matr{Q}_I^{n}$
to be a rotation matrix, and at the same time be explicit in
time. However, the nonlinearity of the constraint on $\matr{Q}_I^{n}$
needs to be solved with an iterative algorithm, which will be addressed
in section \ref{subsection:Res}.

\subsection{Implementation with forces and torques}

For effective implementation of the RATTLE scheme, a difficulty arises
from the fact that we do not necessarily have a 
direct access to $\frac{\partial U_h}{\partial \vect{X}_I}(\vect{X}^{n},
\matr{Q}^{n})$ and $\frac{\partial U_h}{\partial \matr{Q}_I}(\vect{X}^{n}, \matr{Q}^{n})$, as we
compute the expression of forces and torques rather than the functional
$U_h$. In the particular case studied here, we could impose directly
$U_h$ in the computation of velocity and position, but in that case, we
would not be able to treat non-conservative exterior forces and torques, and the
extension of the method to more complex behavior 
laws for the material would become unfeasible. To that end, we have
chosen to recover $\frac{\partial U_h}{\partial \vect{X}_I}(\vect{X}^{n},
\matr{Q}^{n})$ and $\frac{\partial U_h}{\partial
  \matr{Q}_I}(\vect{X}^{n}, \matr{Q}^{n}) $ from the expression of
forces and torques. We prove, in Appendix \ref{sec:implementation}, that
the equations to be solved have the same form as
(\ref{eqn:translation1}--\ref{eqn:rotation3}), replacing $\frac{\partial
  U_h}{\partial \vect{X}_I}$ with $-\vect{\mathcal{F}}_I^n =-
\sum_{J\in\mathcal{V}_I}{\FIJ}$ and $\frac{\partial U_h}{\partial
  \matr{Q}_I}$ with
$-\frac{1}{2}\matr{j}(\vect{\mathcal{M}}_I^n)\matr{Q}_I^n$, where
$\vect{\mathcal{M}}_I^n = \sum_{J\in\mathcal{V}_I}{\MIJ}$, and changing
the Lagrange multipliers.

In order to implement the scheme, without having to compute matrices
$\matr{\Lambda}_I^n$ and $\matr{\tilde{\Lambda}}_I^n$, we follow once more
\cite[Sec VII.5]{Hairer}. We set~:
\begin{align*}
\matr{Y}_I^n &= \transpose{\matr{Q}_I^{n}}
\cdot \matr{P}_I^{n} \\
\matr{Z}_I^{n+1/2} &= \transpose{\matr{Q}_I^{n}} \cdot
\matr{P}_I^{n+1/2} \cdot\matr{D}_I^{-1}
\end{align*}
We use the following algorithm~:
\begin{itemize}
\item We start the time step knowing $\vect{X}_I^{n}$, $\matr{Q}_I^{n}$,
  $\matr{Z}_I^{n-1/2}$ and $\vect{T}_I^{n-1/2}$ (in the first step,
  these last two elements are the null matrix and the null vector).
\item We compute the forces and torques in a submodule of the code,
 using only positions $\vect{X}^{n}$ and $\matr{Q}^{n}$. 
\item The displacement scheme is written~:
\begin{align*}
\vect{T}_I^{n+1/2} &= \vect{T}_I^{n-1/2} +
\Delta t\vect{\mathcal{F}}_I^n \\
\vect{X}_I^{n+1} &= \vect{X}_I^{n} + \frac{\Delta t}{m_I}\vect{T}_I^{n+1/2}
\end{align*}
\item Then, we use the rotation scheme~:
\begin{itemize}
\item Compute $\matr{A}_I^n = \matr{D}_I\cdot\matr{Z}_I^{n-1/2} -
  \transpose{\matr{Z}_I^{n-1/2}}\cdot\matr{D}_I + \Delta
  t\transpose{\matr{Q}_I^{n}}\cdot\matr{j}(\vect{\mathcal{M}}_I^n)\cdot\matr{Q}_I^{n}$ 
\item Find $\matr{Z}_I^{n+1/2}$ such that~:
\begin{equation}
  \left\{
    \begin{array}{l}
      \matr{Id}+\Delta t\matr{Z}_I^{n+1/2} \text{ is orthogonal} \\
      \matr{Z}_I^{n+1/2}\cdot\matr{D}_I-\matr{D}_I\cdot\transpose{\matr{Z}_I^{n+1/2}}
      = \matr{A}_I^n
    \end{array}
  \right.
  \label{eqn:solve1}
\end{equation} 
\item Compute $\matr{Q}_I^{n+1} = \matr{Q}_I^{n}\cdot(\matr{Id}+\Delta t\matr{Z}_I^{n+1/2})$
\end{itemize}
\end{itemize}

We can observe that all those steps are explicit, and that the only
step that requires an iterative resolution is
(\ref{eqn:solve1}). Following \cite{Hairer}, we use the 
quaternion iterative method to solve (\ref{eqn:solve1}) for
$\matr{Z}_{n+1/2}$. We describe that method in the next subsection.

\subsection{Resolution of the nonlinear step} \label{subsection:Res}

Note that $\matr{A}_I^n$ is a skew-symmetric matrix, which can be written
as~:
\begin{equation*}
  \matr{A}_I^n = \left(
   \begin{array}{ccc}
     0 & -\alpha_3 & \alpha_2 \\
     \alpha_3 & 0 & -\alpha_1 \\
     -\alpha_2 & \alpha_1 & 0
   \end{array}
   \right)
\end{equation*}
Equation (\ref{eqn:solve1}) now reads~:
\begin{equation}
  \left\{
    \begin{array}{l}
      \matr{Z}_I^{n+1/2}\cdot\matr{D}_I-\matr{D}_I\cdot\transpose{\matr{Z}_I^{n+1/2}} = \matr{A}_I^n \\
      \left(\matr{Id} +\Delta t\matr{Z}_I^{n+1/2}\right)\cdot\left(\matr{Id}
  +\Delta t\transpose{\matr{Z}_I^{n+1/2}}\right) = \matr{Id} \label{eqn:orthogonal}
    \end{array}
  \right.
\end{equation}
To impose the second line of (\ref{eqn:orthogonal}), we write the matrix
$\matr{Id} +\Delta t\matr{Z}_I^{n+1/2}$ with the quaternion notation~:
\begin{equation*}
  \matr{Id}+\Delta t\matr{Z}_I^{n+1/2} = (e_0^2+e_1^2+e_2^2+e_3^2)\matr{Id} +
  2e_0\matr{E}+2\matr{E}^2 
\end{equation*}
with~:
\begin{equation*}
  E = \left(
  \begin{array}{ccc}
    0 & -e_3 & e_2 \\
    e_3 & 0 & -e_1 \\
    -e_2 & e_1 & 0
  \end{array}
  \right)
\end{equation*}
We make use of the property that every orthogonal matrix can be written in
this form, and that condition $e_0^2+e_1^2+e_2^2+e_3^2 =1$ ensures that such a
matrix is orthogonal. Equation (\ref{eqn:solve1}) is hence equivalent to
solving for $e_0,e_1,e_2,e_3$ the following quadratic system of equations~:
\begin{equation}
  \left\{
  \begin{array}{rl}
    2(d_2+d_3)e_0e_1+2(d_2-d_3)e_2e_3 &= \Delta t\alpha_1 \\
    2(d_1+d_3)e_0e_2+2(d_3-d_1)e_1e_3 &= \Delta t\alpha_2 \\
    2(d_1+d_2)e_0e_3+2(d_1-d_2)e_1e_2 &= \Delta t\alpha_3 \\
    e_0^2+e_1^2+e_2^2+e_3^2 &= 1
  \end{array}
  \right. \label{eqn:quaternions}
\end{equation}

Existence and uniqueness do not hold for this set of equations. In the
simple case where $\alpha_1=\alpha_2=\alpha_3=0$, there are distinct
solutions for $(e_0,e_1,e_2,e_3)$~: $(1,0,0,0)$ (in that case,
$\matr{Z}^{n+\frac{1}{2}}=\matr{Id}$), $(0,1,0,0)$ (in that case,
$\matr{Z}^{n+\frac{1}{2}}$ represents the axial symmetry around axis
$x$), $(0,0,1,0)$ (associated with the axial symmetry around axis
$y$), $(0,0,0,1)$ (associated with the axial symmetry around axis
$z$), and their opposites which represent the same
transformation. There is a deep physical reason for that
non-uniqueness~: dynamically speaking, the rigid body is totally
represented by its equivalent inertia ellipsoid (the ellipsoid with
the same axes of inertia and moments of inertia), which 
is invariant under the axial symmetries around the inertial axes $x$,
$y$ and $z$. As the rotation $\matr{Id}+\Delta t\matr{Z}_I^{n+1/2}$ is
an increment of the global rotation of the particle, we select a
solution ``close'' to identity, in a certain sense.

The existence and uniqueness in a neighbourhood of identity can be
obtained from the equivalent formulation of RATTLE using the discrete
Moser-Veselov scheme, with a fixed point theorem applied on equation (17)
of reference \cite{Vilmart}. We have found an explicit bound on
the time-step $\Delta t$ for the iterative scheme to
converge, and ensure existence and uniqueness in a neighbourhood of
identity. It is derived in Appendix \ref{sec:existence}. We use the following iterative scheme \cite{Hairer}~:
\begin{itemize}
\item We start with $(e_0^0,e_1^0,e_2^0,e_3^0) = (1,0,0,0)$ (which
  represents identity).
\item At each iteration, we compute~:
\begin{align}
  e_1^{k+1} &= \frac{\Delta
    t\alpha_1-2(d_2-d_3)e_2^ke_3^k}{2(d_2+d_3)e_0^k} \label{eqn:iterative1} \\
  e_2^{k+1} &= \frac{\Delta
    t\alpha_2-2(d_3-d_1)e_1^ke_3^k}{2(d_1+d_3)e_0^k} \label{eqn:iterative2} \\
  e_3^{k+1} &= \frac{\Delta
    t\alpha_3-2(d_1-d_2)e_1^ke_2^k}{2(d_1+d_2)e_0^k} \label{eqn:iterative3} \\
  e_0^{k+1} &= \sqrt{1-(e_1^{k+1})^2-(e_2^{k+1})^2-(e_3^{k+1})^2}
  \label{eqn:iterative4} 
\end{align}
\end{itemize}
Let us introduce~:
\begin{equation*}
  \mathcal{B}(\frac{\sqrt{2}}{2}) = \left\{(e_0,e_1,e_2,e_3)/e_0^2+e_1^2+e_2^2+e_3^2=1,e_1^2+e_2^2+e_3^2<\frac{1}{2} \right\}
\end{equation*}
When the time-step $\Delta t$ satisfies the condition~:
\begin{equation}
  \Delta
  t\left(\frac{|\alpha_1|}{I_1}+\frac{|\alpha_2|}{I_2}+\frac{|\alpha_3|}{I_3}\right)\leq\frac{\sqrt{21}-3}{6}\approx
  0.26 \label{eqn:cond_CFL}
\end{equation}
the algorithm (\ref{eqn:iterative1})--(\ref{eqn:iterative4}) converges
with a geometrical speed 
to the unique solution in $\mathcal{B}(\frac{\sqrt{2}}{2})$.

Let us observe that $I_i$ and $\matr{D}$ scale as $\rho h^5$. In
addition, as $\matr{P}_I =
\matr{j}(\vect{\Omega}_I)\matr{Q}_I\matr{D}_I$,
$\matr{Z}_I^{n+\frac{1}{2}}$ is of the order of
$\lVert\vect{\Omega}_I\rVert$. Using the expressions (\ref{eqn:torque}) and
(\ref{eqn:bendingtorsion}), and the fact that $\alpha_n$, $\alpha_s$
and $\alpha_t$ scale as $h^2$, we obtain that $\vect{\mathcal{M}}_I^n$
is of the order of $E h^3$. Condition (\ref{eqn:cond_CFL}) therefore
gives us a constraint on the time-step of the following type~:
\begin{equation}
  \Delta t \lVert\vect{\Omega}_I\rVert + \frac{\Delta
  t^2}{h^2}\frac{E}{\rho}\leq C
\end{equation} 
where $C$ is a constant. This is the natural CFL condition for an
explicit scheme on
rotation, with $\sqrt{\frac{E}{\rho}}$ the typical celerity of the
compression and shear waves in the material.

\section{Numerical results}

In this section, we present several challenging test cases. First, we
address Lamb's problem, which allows us to examine numerically the
precision of the method in the case of small displacements against a
semi-analytic solution. The presence of surface waves is the most
difficult part of the problem, and the results appear to be satisfactory. 
We examine the conservation of energy on
the case of a three-dimensional cylinder submitted to large
displacement. In the end, we also demonstrate the ability of the method
to tackle static rod and shell problems using the same formulation, on the cases
of the bending of a rod and of the loading of a hemispherical shell.

\label{sec:resultsnum}

\subsection{Lamb's problem}

We have simulated Lamb's problem (see \cite{Lamb})~: a
semi-infinite plane is described by a rectangular domain,
with a free surface on the upper side, and absorbing conditions on the
other sides. On a surface particle, we apply a vertical force, whose
time evolution is described by a Ricker function (the second
derivative of a Gaussian function). We observe the propagation of
three waves~: inside the domain, a compression wave of type P and a
shear wave of type S, and on the surface, a Rayleigh wave. We also
have a P-S wave linking the P and the S waves, which is a conversion
of the P wave into an S wave after reflection at the surface. In the case of a
two-dimensional problem, the intensity of P and S waves is inversely
proportional to the distance to the source, and the intensity of the
Rayleigh wave is preserved throughout its propagation.


We have chosen the following characteristics for the material~: the
density is $\rho = 2200\text{ kg.m}^{-3}$, the Poisson coefficient is
$\nu = 0.25$, Young's modulus is $E = 1.88.10^{10}\text{ Pa}$. The
velocity of P waves is therefore approximately $3202\text{ m.s}^{-1}$
and the velocity of S waves is $1849\text{ m.s}^{-1}$.

The force applied is a Ricker of
central frequency $14.5\text{ Hz}$, that is, with maximal frequency
around $40\text{ Hz}$. The minimal wave length for P waves is
therefore $80\text{ m}$, and the minimal wave length for S waves is
approximately $50\text{ m}$. In the rest of this subsection, we
call ``wave length'' this minimal wave length of  $50\text{ m}$. We
indicate the discretization step in terms of number of elements
per wave length.

\begin{figure}[h!]
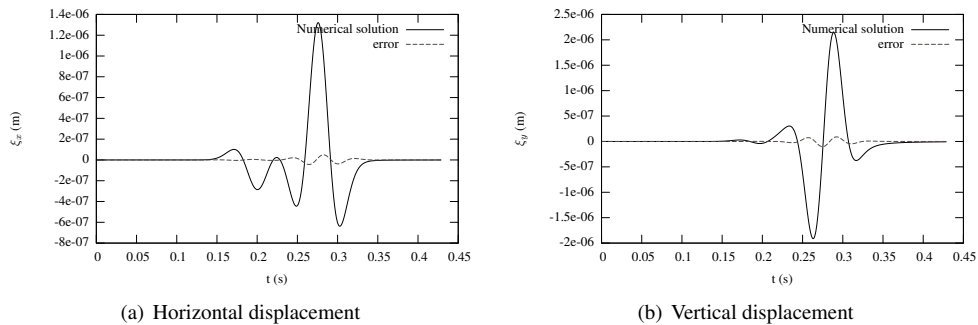

  \centering
  \subfigure[Horizontal displacement]{
    \resizebox{0.45\textwidth}{!}{
      \input{Ux155mod.tex}
    }
  \label{fig:Ux155} }
  \subfigure[Vertical displacement]{
    \resizebox{0.45\textwidth}{!}{
      \input{Uy155mod.tex}
    }
  \label{fig:Uy155} }
  \caption{Displacement at the surface, 300 meters from source, with $\Delta x=5\text{ m}$,
  $\Delta y=5\text{ m}$ (10 points per wave length)} \label{fig:U155}
\end{figure}

Lamb's problem has the interesting particularity of having a
semi-analytic solution~: Cagniard's method is described in
\cite{DeHoop}. We have compared our results with 
this exact solution and thus estimate the numerical error of the
scheme. The comparison between the numerical results and the
semi-analytic solution obtained at 300 meters
from the source, on the surface, with $\Delta x = \Delta y = 5\text{
  m}$ (10 points per wave length), is shown on figure \ref{fig:U155}. 

We compute the same result with different spatial discretizations, with
$\Delta x = \Delta y$. As expected, refining the spatial discretization
decreases the error. The velocity of the different waves agrees with
the exact solution, and the amplitude of the waves is accurately
captured with more than 10 elements per wavelength.
The accuracy of the method cannot compare with that of spectral elements (5
points per wave length), but it gives better results than classic
second-order finite elements (30 points per wave length), and mostly
on the surface, where we recover the non-dissipative Rayleigh wave. This is
probably due to the introduction of parameter $\vect{\theta}$ which helps us
simulate the rotation of the particle precisely, instead of recovering
it as a Taylor development of the displacement, thus losing one
order of accuracy for rotation.

\begin{figure}[h!]
  \centering
  \resizebox{0.65\textwidth}{!}{
    \input{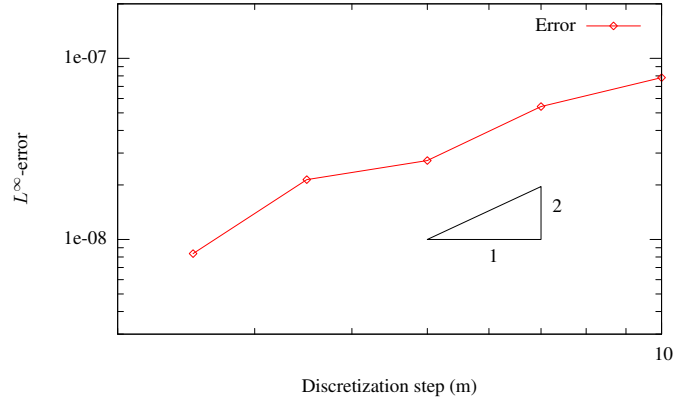}
  }
  \caption{Linear fitting of the log-log diagram for the numerical
  error against the spatial discretization step} \label{fig:logerreur}
\end{figure}

If we measure the $L^\infty$-error on vertical displacement at 300
meters from the source, with an angle of 60\textdegree~with the
horizontal axis, we obtain an approximate slope of 2 fitting the
points (figure \ref{fig:logerreur}). This confirms the results of
subsection \ref{sec:equivalent} as to the second-order nature of the
spatial scheme. 

\subsection{Conservation of energy}

\begin{figure}[h!]
  \centering
  \includegraphics[width=7cm,angle=270]{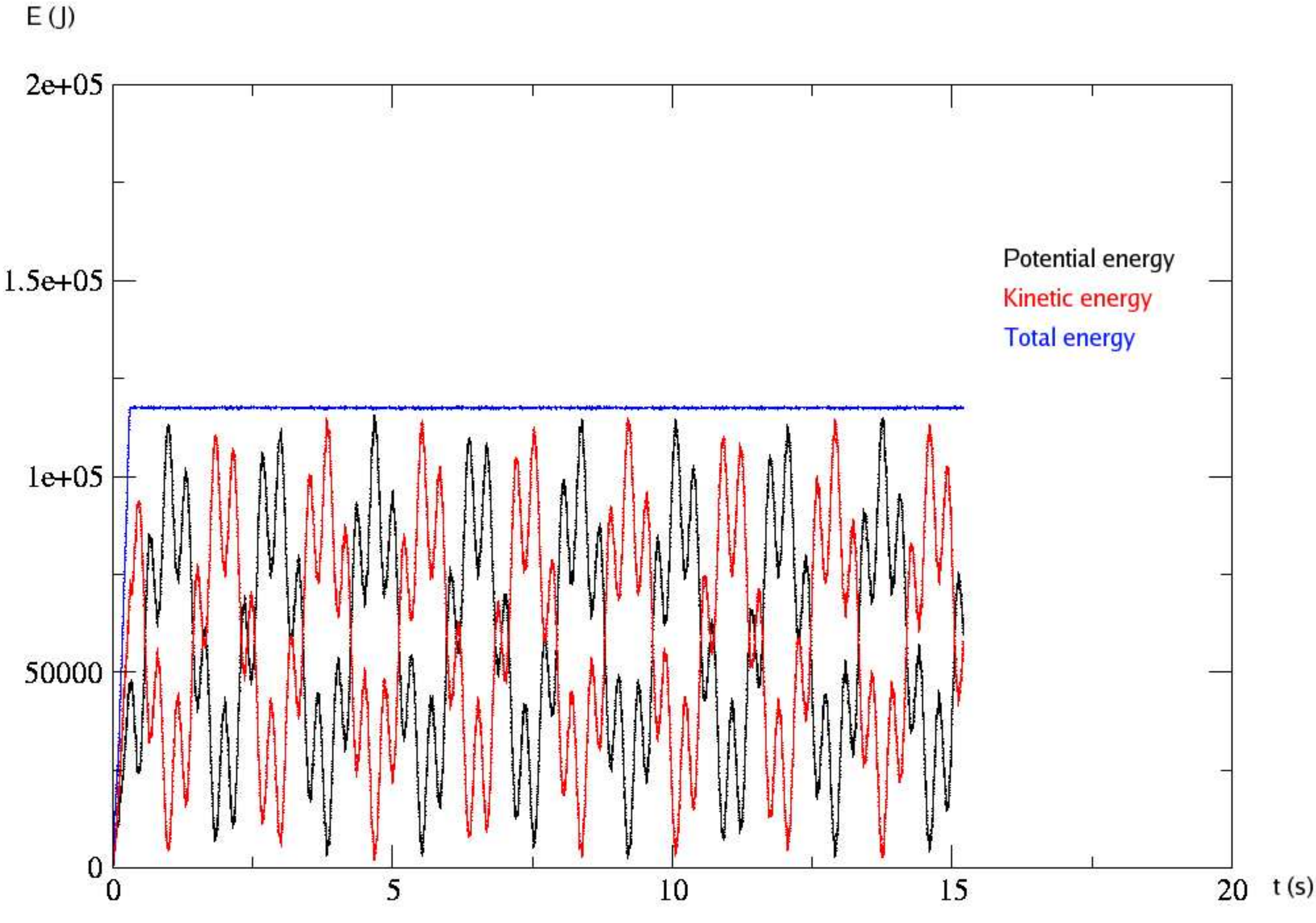}
  \caption{Total, potential and kinetic energies for the simulation of
    the cylinder over 500000 time-steps} \label{fig:energycylinder}
\end{figure}

\begin{figure}[h!]
  \centering
  \includegraphics[height=7cm]{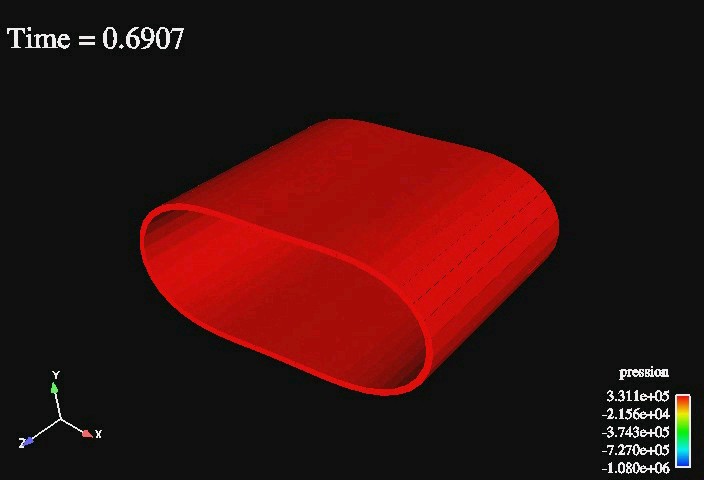}
  \caption{Initial configuration of the cylinder}
  \label{fig:cylinder}
\end{figure}

In order to illustrate the conservation of energy by the scheme, we model the
evolution of a pinched cylinder. The cylinder has a radius of 1m, a
height of 2m and a width of 1cm. The physical characteristics are that
of steel ($E=210000\text{ MPa}$, $\nu=0.25$). The cylinder is
discretized with 50 elements on the perimeter, 20 elements on the
height and one element in width. Opposite forces are applied on two
sides of the cylinder, pinching it. At the initial time, the forces
are removed, and the cylinder is left free. We simulate the system
over 500,000 time-steps, corresponding to 45 oscillations of the first
mode of the cylinder. The large number of time-steps required reflects
the fact that a number of smaller local oscillations propagate at high
velocities, and that the cylinder is very thin. On figure
\ref{fig:energycylinder}, we observe an excellent preservation of the
energy. The configuration of the cylinder at the moment of release is
shown on figure \ref{fig:cylinder}.

The preservation of energy is quite satisfactory, even with large
displacements in a three-dimensional geometry.

\subsection{Static shell test cases}

In order to show the versatility of the method, we compare the static
deformation obtained with Mka3D (adding damping to the model) to the
second and fourth benchmarks for geometric nonlinear shells found in
\cite{Sze}. 

The first benchmark considered is that of the cantilever subjected to an
end moment $\mathcal{M}$. Let $N$ be the number of discrete elements in
the length of the beam. We take one element in the two other
directions. We immediately see that at the equilibrium, for each
particle $I$, the sum of forces is null, and using the boundary
conditions, the force $\FIJ$ between particles is always null. The sum
of moments is also zero, and is equal to the end moment
$\mathcal{M}$. As $\FIJ=\vect{0}$, if we denote $\theta_N$ the angle
between two consecutive particles, using (\ref{eqn:alpha}),
\begin{equation}
  \MIJ = \MIJ^f = \frac{EI}{2\DIJeq}\sin\theta_N
\end{equation} 
If we take the maximum end moment $\mathcal{M}_{max}=2\pi\frac{EI}{L}$,
which is the theoretical moment applied to bend the beam into a circle,
we obtain~:
\begin{equation}
  N\theta_N = N\arcsin\left(\frac{2\pi}{N}\right)
\end{equation}
As $N$ tends to infinity, the deflection angle of the end $N\theta_N$
tends to $2\pi$ with second order precision, which indicates a second
order convergence to the theoretical solution. This convergence has been
checked in practice.

\begin{figure}[h!]
  \centering
  \includegraphics[width=7cm]{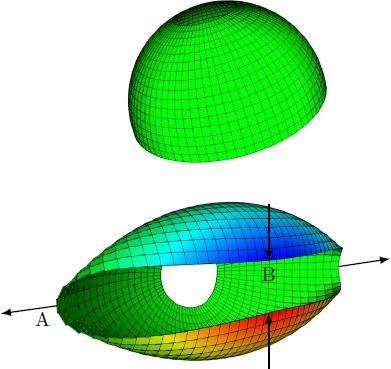}
  \caption{Initial geometry and deformed geometry at $\mathcal{F}=400N$
    for the hemispherical shell subjected to alternating radial forces}
  \label{fig:sphere}
\end{figure}

\begin{figure}[h!]
  \centering
  \includegraphics[width=7cm]{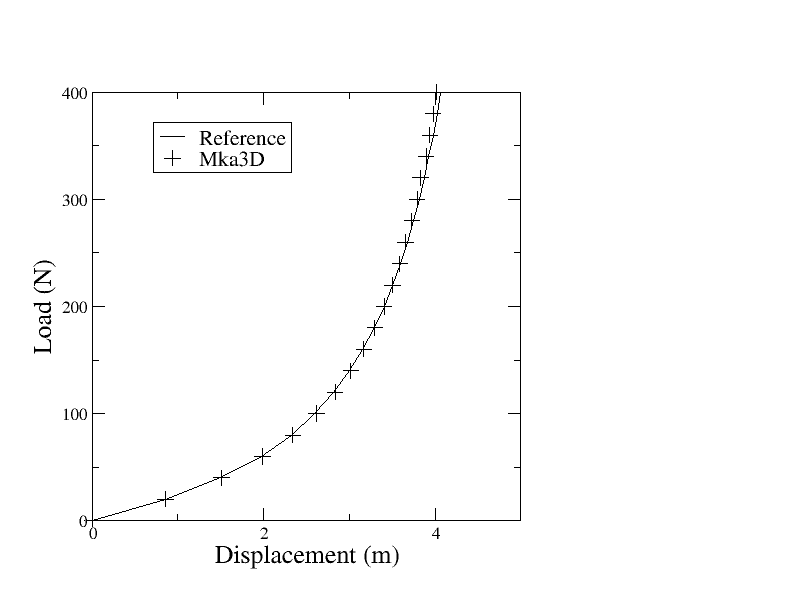}
  \includegraphics[width=7cm]{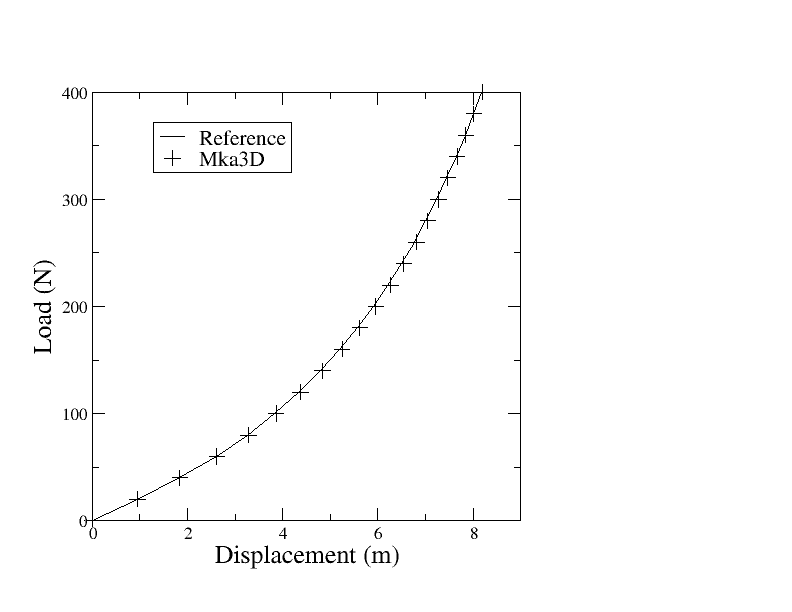}
  \caption{Load-deflection curves for the hemisphere shell at the
    loading points A (left) and B (right)}
  \label{fig:load}
\end{figure}

The second benchmark considered is a hemispherical shell with an
$18^{\circ}$ circular cutout at its pole, loaded by alternating radial point
forces $\mathcal{F}$ at $90^{\circ}$ intervals. The shell is discretized by
16 elements in latitude, 64 elements in longitude and one element in
thickness. The initial and deformed geometries are shown on figure
\ref{fig:sphere}. The radial deflections at the points of loading A and B are
compared with the results obtained in \cite{Sze} in figure
\ref{fig:load}. Our results are in very good agreement with the benchmark.

\section{Conclusion}
\label{sec:conclusion}

In this paper, we proposed a numerical discretization of material
continuum, allowing for the simulation of three-dimensional wave
propagation as well as shell or multibody dynamics, in a monolithic way. It is
consistent with the equations of elastodynamics at order 2 in space
and in time, and we numerically recover the propagation of seismic
waves in the body of the material and at the free surface. Furthermore,
the dynamics of the system are written in the form of a Hamiltonian
dynamics. Using symplectic schemes, we correctly reproduce the 
preservation of the system energy. This ensures
numerical $L^2$-stability of the scheme, and allows long-time stable
simulations with large displacements and large deformations. As the
method is entirely local and requires no matrix inversion, it can be
easily parallelized with domain decomposition. The main restriction is
the size of the time-step due to the explicit nature of the
integration scheme. This could be remedied by using asynchronous
symplectic integrators in order to have local time refinement at small
elements and a global larger time-step \cite{LewMarsden}. This work
can be seen as a first step towards using more complex constitutive laws
(while still maintaining stability of the scheme), and towards coupling
particle dynamics simulation with a fluid dynamics simulation for
fluid-structure interaction.

\acks{
The first author acknowledges the support of CEA under Grant
n\ensuremath{^\circ}1045. 

We would like to thank Serge Piperno, Tony Leli\`evre, Fr\'ed\'eric
Legoll and Eric Canc\`es (Cermics and UR Navier, Ecole des Ponts) for useful
discussions and advice on the mathematical and computational aspects of
this paper. We also thank Karam Sab (UR Navier, Ecole des Ponts) for
pointing us the similarity of our model with Cosserat 
models. Thanks are also due to Gilles Vilmart for discussions on 
the resolution of the quaternion scheme.}

\appendix
\setcounter{section}{0}

\section{Expression of the equivalent volumetric deformation with a free
surface}
\label{sec:freesurf}

We need to account for the boundary condition
$\matr{\sigma}\cdot\vect{n}=\vect{0}$ at every free surface of the
particles. We have seen in section \ref{sec:derivation} that the discrete equivalent for
$\matr{\sigma}\cdot\vect{n}$ is $\FIJ$. For a given particle $I$, we
assume that the particle is surrounded by real particles
$J\in\mathcal{V}_I$, and by `ghost' particles $J\in\mathcal{V}_I^l$ at
every free boundary. The position of these particles is ajusted in order
to satisfy the boundary condition.

The equivalent deformation of particle $I$ can be expressed as in the
bulk of the material~:
\begin{equation*}
  \varepsilon_I^v =
  \sum_{J\in\mathcal{V}_I}{\frac{1}{2}\frac{\SIJ}{V_I}\DeltauIJ\cdot\nIJ}
  + \sum_{J\in\mathcal{V}_I^l}{\frac{1}{2}\frac{\SIJ}{V_I}\DeltauIJ\cdot\nIJ}
\end{equation*}
For a ghost particle $J\in\mathcal{V}_I^l$, the boundary condition
$\FIJ\cdot\nIJ=0$ boils down to~:
\begin{equation}
  \frac{\SIJ}{\DIJeq}\frac{E}{1+\nu}\DeltauIJ\cdot\nIJ+\SIJ\frac{E\nu}{(1+\nu)(1-2\nu)}\varepsilon_I^v=0
\label{eqn:BC}
\end{equation}
Summing (\ref{eqn:BC}) over the ghost particles, and using the fact
that the free volume $V_I^l$ satisfies
\begin{equation*}
  V_I^l = \sum_{J\in\mathcal{V}_I^l}{\frac{\SIJ\DIJeq}{6}}
\end{equation*}
we find that the deformation of the links with the ghost particles
should follow the equation~:
\begin{equation*}
  \sum_{J\in\mathcal{V}_I^l}{\frac{1}{2}\frac{\SIJ}{V_I}\DeltauIJ\cdot\nIJ}
  = -\frac{3\nu}{1-2\nu}\frac{V_I^l}{V_I+\frac{3\nu}{1-2\nu}V_I^l}\sum_{J\in\mathcal{V}_I}{\frac{1}{2}\frac{\SIJ}{V_I}\DeltauIJ\cdot\nIJ}
\end{equation*}
Inserting this relation in the expression of $\varepsilon_I^v$, we check
that~:
\begin{equation*}
  \varepsilon_I^v =
  \sum_{J\in\mathcal{V}_I}{\frac{1}{2}\frac{\SIJ}{V_I+\frac{3\nu}{1-2\nu}V_I^l}\DeltauIJ\cdot\nIJ}
\end{equation*}

\section{Expression of the coefficients for the flexion and torsion of
  the particle links}
\label{sec:flexion}

We denote~:
\begin{align}
  I_{IJ}^s &= \iint_{\mathcal{S}_{IJ}}{(\vect{X\PIJ}\cdot\sIJ)^2 dX}
  \\
  I_{IJ}^t &= \iint_{\mathcal{S}_{IJ}}{(\vect{X\PIJ}\cdot\tIJ)^2 dX}
\end{align} 
the principal moments of the interface between particles $I$ and $J$,
we require that~:
\begin{equation}
  \left\{
  \begin{array}{l}
    \displaystyle \alpha_n+\alpha_s = \frac{EI_{IJ}^s}{\SIJ} \label{eqn:alpha}\\
    \displaystyle \alpha_n+\alpha_t = \frac{EI_{IJ}^t}{\SIJ} \\
    \displaystyle \alpha_s+\alpha_t = \frac{E(I_{IJ}^s+I_{IJ}^t)}{2(1+\nu)\SIJ} 
  \end{array}
  \right.
\end{equation}
The expression of the $\alpha$ is given by~:
\begin{align}
  \alpha_n &= \frac{(1+2\nu)E}{4(1+\nu)\SIJ}(I_{IJ}^s+I_{IJ}^t)\\
  \alpha_s &=
  \frac{E}{4(1+\nu)\SIJ}((3+2\nu)I_{IJ}^s-(1+2\nu)I_{IJ}^t)\\
  \alpha_t &=
  \frac{E}{4(1+\nu)\SIJ}((3+2\nu)I_{IJ}^t-(1+2\nu)I_{IJ}^s)
\end{align}

\section{Derivation of the forces and torques from the potential
  energies}
\label{sec:forces}

The derivation of potential energies is straightforward~:
\begin{align*}
  \frac{\partial U_t}{\partial\vect{X}_I} &=
  -\sum_{J\in\mathcal{V}_I}{\frac{\SIJ}{\DIJeq}\frac{E}{1+\nu}\DeltauIJ}
  \\
  \frac{\partial U_d}{\partial\vect{X}_I} &=
  -\sum_{J\in\mathcal{V}_I}{\frac{E\nu}{(1+\nu)(1-2\nu)}\SIJ\epsilonIJ\left(\nIJ+\frac{1}{\DIJ}\DeltauIJ-\frac{1}{\DIJ}(\DeltauIJ\cdot\nIJ)\nIJ\right)}
  \\
  \frac{\partial U_t}{\partial \matr{Q}_I} &=
  -\sum_{J\in\mathcal{V}_I}{\frac{\SIJ}{\DIJeq}\frac{E}{1+\nu}\DeltauIJ\otimes\vect{X_I^0\PIJ}}
  \\
  \frac{\partial U_d}{\partial \matr{Q}_I} 
  &=
  -\sum_{J\in\mathcal{V}_I}{\frac{E\nu}{(1+\nu)(1-2\nu)}\SIJ\epsilonIJ\nIJ\otimes\vect{X_I^0\PIJ}}
  \\
  \frac{\partial U_f}{\partial \matr{Q}_I} &=
  -\sum_{J\in\mathcal{V}_I}{\SIJ\frac{E}{\DIJeq}\left(\alpha_n(\matr{Q}_J\cdot\nIJ^0)\otimes\nIJ^0 + \alpha_s(\matr{Q}_J\cdot\sIJ)\otimes\sIJ + \alpha_t(\matr{Q}_J\cdot\tIJ)\otimes\tIJ \right)}
\end{align*}

Using the expression of the force $\FIJ$ between particles $I$ and $J$~:
\begin{equation*}
  \FIJ = \frac{\SIJ}{\DIJeq}\frac{E}{1+\nu}\DeltauIJ 
  + \SIJ\frac{E\nu}{(1+\nu)(1-2\nu)}\epsilonIJ\left(\nIJ+\frac{1}{\DIJ}\DeltauIJ-\frac{1}{\DIJ}(\DeltauIJ\cdot\nIJ)\nIJ\right)
\end{equation*}
we obtain~:
\begin{equation*}
  m_I\vect{\dot{v}}_I = \vect{\dot{T}}_I = \FIJ
\end{equation*}

For the rotational part, it is easily obtained that~:
\begin{equation*}
  \matr{j}(\matr{R}_I\vect{\Omega}_I) =
  \matr{j}(\vect{\Omega}_I)\matr{D} - \matr{D}\matr{j}(\vect{\Omega}_I)
  = \matr{P}_I\cdot\transpose{\matr{Q}} - \matr{Q}_I\cdot\transpose{\matr{P}}
\end{equation*}
Deriving in time, we obtain~:
\begin{equation*}
  \frac{d}{dt}\left(\matr{j}(\matr{R}_I\cdot\vect{\Omega}_I)\right) =
  -\left(\frac{\partial
      H_h}{\partial\matr{Q}_I}\right)\transpose{\matr{Q}_I} + \matr{Q}_I\transpose{\left(\frac{\partial H_h}{\partial\matr{Q}_I}\right)}
\end{equation*}
Using the fact that~:
\begin{equation*}
  (\vect{a}\otimes\vect{b})\cdot\matr{Q} = \vect{a}\otimes(\transpose{\matr{Q}}\cdot\vect{b})
\end{equation*}
we get~:
\begin{equation}
  \frac{\partial U_t}{\partial
  \matr{Q}_I}\cdot\transpose{\matr{Q}_I} 
  =
  -\sum_{J\in\mathcal{V}_I}{\frac{\SIJ}{\DIJeq}\frac{E}{1+\nu}\DeltauIJ\otimes(\matr{Q}_I\cdot\vect{X_I^0\PIJ})} \label{eqn:rotationt}
\end{equation}
\begin{equation}
  \frac{\partial U_d}{\partial
  \matr{Q}_I}\cdot\transpose{\matr{Q}_I} 
  = -\sum_{J\in\mathcal{V}_I}{\frac{E\nu}{(1+\nu)(1-2\nu)}\SIJ\epsilonIJ\nIJ\otimes(\matr{Q}_I\cdot\vect{X_I^0\PIJ})}\label{eqn:rotationd}
\end{equation}
\begin{multline}
  \frac{\partial U_f}{\partial
  \matr{Q}_I}\cdot\transpose{\matr{Q}_I}
  =
  -\sum_{J\in\mathcal{V}_I}{\SIJ\frac{E}{\DIJeq}\left(\alpha_n(\matr{Q}_J\cdot\nIJ^0)\otimes(\matr{Q}_I\cdot\nIJ^0)
      +
      \alpha_s(\matr{Q}_J\cdot\sIJ)\otimes(\matr{Q}_I\cdot\sIJ)\right.}\\ 
  {\left.+
      \alpha_t(\matr{Q}_J\cdot\tIJ)\otimes(\matr{Q}_I\cdot\tIJ)\right)} \label{eqn:rotationflex}
\end{multline}

Denoting $\text{symm}()$ and $\text{skew}()$ the symmetric and
skew-symmetric parts of a matrix, we note that for any $\vect{a}$ and
$\vect{b}$~:
\begin{equation*}
  \matr{j}(\vect{a}\wedge\vect{b}) = -\text{skew}(\vect{a}\otimes\vect{b})
\end{equation*}
Using the expression of the torques $\MIJ^t$ and $\MIJ^f$~:
\begin{equation*}
  \MIJ^t = \frac{\SIJ}{\DIJeq}\frac{E}{1+\nu}
  (\matr{Q}_I\cdot\vect{X_I^0\PIJ})\wedge\DeltauIJ
  + \frac{E\nu}{(1+\nu)(1-2\nu)}\epsilonIJ\SIJ(\matr{Q}_I\cdot\vect{X_I^0\PIJ})\wedge\nIJ
\end{equation*}
\begin{equation*}
  \MIJ^f = \frac{\SIJ}{\DIJeq}\left(\alpha_n
  (\matr{Q}_I\cdot\nIJ^0)\wedge(\matr{Q}_J\cdot\nIJ^0)\right. 
  \left. + \alpha_s (\matr{Q}_I\cdot\sIJ)\wedge(\matr{Q}_J\cdot\sIJ)
  \right. 
  \left. + \alpha_t (\matr{Q}_I\cdot\tIJ)\wedge(\matr{Q}_J\cdot\tIJ)
  \right)
\end{equation*}
equation (\ref{eqn:rotvit}) gives us the equation on the
angular velocity~:
\begin{equation*}
  \frac{d}{dt}\left(\matr{R}_I\cdot\vect{\Omega}_I\right) =
\sum_{J\in\mathcal{V}_I}{\MIJ^t+\MIJ^f}
\end{equation*}

\section{Details on the implementation of the RATTLE scheme with forces
  and torques}
\label{sec:implementation}

For forces, the relation is simple~:
\begin{equation*}
  \frac{\partial U_h}{\partial \vect{X}_I}(\vect{X},\matr{Q}) = -\sum_{J\in\mathcal{V}_I}{\FIJ}
\end{equation*}

For torques, we have~:
\begin{equation*}
\frac{\partial U_h}{\partial \matr{Q}_I}(\vect{X}, \matr{Q}) = \dot{\matr{P}}_I -
\matr{Q}_I \matr{\Lambda}_I
\end{equation*}
where $\matr{\Lambda}_I$ is the symmetric matrix of Lagrange multipliers
associated with constraint $\matr{Q}_I\cdot\transpose{\matr{Q}_I}
= \matr{Id}$. On the other hand, 
\begin{align*}
  \matr{j}\left(\sum_{J\in\mathcal{V}_I}{\vect{M}_{IJ}}\right)
  = & \matr{\dot{P}}_I\cdot\transpose{\matr{Q}_I}
  +\matr{P}_I\cdot\transpose{\matr{\dot{Q}}_I} 
  -\matr{\dot{Q}}_I\cdot\transpose{\matr{P}_I}
  - \matr{Q}_I\cdot\transpose{\matr{\dot{P}}_I}
  \\
  = & \matr{Q}_I\cdot \transpose{\left(\frac{\partial U_h}{\partial \matr{Q}_I}(\vect{X}, \matr{Q})\right)} 
  - \frac{\partial U_h}{\partial \matr{Q}_I}(\vect{X}, \matr{Q})\cdot\transpose{\matr{Q}_I}
\end{align*}
as the $\matr{\Lambda}_I$ are symmetric. Therefore,
there exists a symmetric matrix $\matr{\Lambda}^0_I$ such that~:
\begin{equation*}
  \frac{\partial U_h}{\partial \matr{Q}_I}(\vect{X}, \matr{Q})
  = \left(-\frac{1}{2}\matr{j}\left(\sum_{J\in\mathcal{V}_I}{\vect{M}_{IJ}}\right)-\matr{\Lambda}^0_I\right)\cdot\matr{Q}_I
\end{equation*}

We denote~:
\begin{align*}
\vect{\mathcal{F}}_I^n &=
\sum_{J\in\mathcal{V}_I}{\FIJ} \\
\vect{\mathcal{M}}_I^n &= \sum_{J\in\mathcal{V}_I}{\MIJ}
\end{align*}
where forces $\FIJ$ and torques $\MIJ$ have been
computed with positions $\vect{X}^{n}$ and $\matr{Q}^{n}$.

We can rewrite equations (\ref{eqn:translation1}) to
(\ref{eqn:rotation3}) as follows~:
\begin{align}
  \vect{T}_I^{n+1/2} &= \vect{P}_I^{n} + \frac{\Delta
  t}{2}\vect{\mathcal{F}}_I^n  \label{eqn:translation1modif} \\
  \matr{P}_I^{n+1/2} &= \matr{P}_I^{n} + \frac{\Delta
  t}{4}\matr{j}(\vect{\mathcal{M}}_I^n)\matr{Q}_I^{n} 
  + \frac{\Delta
  t}{2}(\matr{\Lambda}_I^n+\matr{\Lambda}_I^{n,0})\matr{Q}_I^{n} 
  \label{eqn:rotation1modif}
\end{align}
\begin{align}
  \vect{X}_I^{n+1} &= \vect{X}_I^{n} + \frac{\Delta
  t}{m_I}\vect{T}_I^{n+1/2} \label{eqn:translation2modif} \\
  \matr{Q}_I^{n+1} &= \matr{Q}_I^{n} +
  \Delta t\matr{P}_I^{n+1/2}\matr{D}_I^{-1}\label{eqn:rotation2modif}
\end{align}
\begin{equation}
  \text{where } \matr{\Lambda}_I^n \text{ is such that
  }\transpose{\matr{Q}_I^{n+1}}\cdot\matr{Q}_I^{n+1} = \matr{Id}  
\end{equation}
\begin{align}
  \vect{T}_I^{n+1} &= \vect{T}_I^{n+1/2} + \frac{\Delta
  t}{2}\vect{\mathcal{F}}_I^{n+1} \label{eqn:translation3modif}\\
  \matr{P}_I^{n+1} &= \matr{P}_I^{n+1/2} + \frac{\Delta
    t}{4}\matr{j}(\vect{\mathcal{M}}_I^{n+1})\matr{Q}_I^{n+1} 
  + \frac{\Delta t}{2}(\matr{\tilde{\Lambda}}_I^{n+1}+\matr{\tilde{\Lambda}}_I^{n+1,0})\matr{Q}_I^{n+1}
  \label{eqn:rotation3modif} ,
\end{align}
\begin{equation}
  \text{where } \matr{\Lambda}_I^n  \text{ is such that }\transpose{\matr{Q}_I^{n+1}}\cdot\matr{P}_I^{n+1}\cdot\matr{D}_I^{-1} +
\matr{D}_I^{-1}\cdot\transpose{\matr{P}_I^{n+1}}\cdot\matr{Q}_I^{n+1} = \matr{0}
\end{equation}

\section{Resolution of the nonlinear step of the RATTLE
  time-scheme}
\label{sec:existence}

In this appendix, we examine the resolution of the nonlinear step of
the RATTLE time-scheme described in section \ref{subsection:Res}. We
determine conditions on the time-step $\Delta t$ that ensure 
convergence of the iterative algorithm
(\ref{eqn:iterative1})--(\ref{eqn:iterative4}) in a certain
neighbourhood of 
identity, and we conclude on the existence and uniqueness of a solution in
this neighbourhood.

We denote $\mathcal{B}(\vect{0},r)$ the
ball of center $\vect{0}$ and radius $r$~:
\begin{equation*}
  \mathcal{B}(\vect{0},r) = \left\{(e_1,e_2,e_3)/e_1^2+e_2^2+e_3^2<r^2 \right\}
\end{equation*}
Using the numerical scheme described in section \ref{subsection:Res}, we
first show that it stabilizes a ball included in
$\mathcal{B}(\vect{0},\frac{\sqrt{2}}{2})$, under a CFL-type condition
on $\Delta t$. We then show convergence in that same ball, and we
conclude on convergence to the unique fixed point. 

\subsection{The iterative scheme is bounded}

Starting with a given $(e_0,e_1,e_2,e_3)$ computed in the previous
iteration, such that $e_0^2+e_1^2+e_2^2+e_3^2=1$, the iterative scheme
(\ref{eqn:iterative1})--(\ref{eqn:iterative4}) 
gives the new quadruplet $(e_0^*,e_1^*,e_2^*,e_3^*)$ defined by~:
\begin{align*}
  e_1^* &= \frac{\Delta t\alpha_1-2(d_2-d_3)e_2e_3}{2(d_2+d_3)e_0} \\
  e_2^* &= \frac{\Delta t\alpha_2-2(d_3-d_1)e_1e_3}{2(d_1+d_3)e_0} \\
  e_3^* &= \frac{\Delta t\alpha_3-2(d_1-d_2)e_1e_2}{2(d_1+d_2)e_0} \\
  e_0^* &= \sqrt{1-(e_1^*)^2-(e_2^*)^2-(e_3^*)^2}
\end{align*}

For this scheme to be well-defined, $(e_1^*,e_2^*,e_3^*)$ should be in
$\mathcal{B}(\vect{0},1)$. We impose a stronger condition, with
$(e_1,e_2,e_3)$ and $(e_1^*,e_2^*,e_3^*)$ in $\mathcal{B}(\vect{0},\beta)$
where $\beta$ is less than $\frac{1}{2}$.

Suppose that~:
\begin{equation*}
  e_1^2+e_2^2+e_3^2<\beta
\end{equation*}
We want to have~:
\begin{equation*}
  (e_1^*)^2+(e_2^*)^2+(e_3^*)^2<\beta
\end{equation*}
As $e_0^2+e_1^2+e_2^2+e_3^2=1$, we also have $e_0^2>1-\beta$. Since~:
\begin{equation*}
  |e_2e_3|\leq\frac{1}{2}(e_2^2+e_3^2)<\frac{\beta}{2}
\end{equation*}
we obtain~:
\begin{equation*}
  |e_1^*|<\frac{1}{2\sqrt{1-\beta}(d_2+d_3)}(|\Delta t\alpha_1|+\beta|d_2-d_3|)
\end{equation*}
Let us define $I_1=d_2+d_3$, $I_2=d_1+d_3$, $I_3=d_1+d_2$ and~:
\begin{multline*}
  f(\beta) =
  \frac{1}{4(1-\beta)}\left[\Delta t^2\left(\frac{|\alpha_1|^2}{I_1^2}
  + \frac{|\alpha_2|^2}{I_2^2}
  + \frac{|\alpha_3|^2}{I_3^2}\right)\right. 
  \left.+ 2\beta \Delta t\left(\frac{|d_2-d_3||\alpha_1|}{I_1^2}
  + \frac{|d_3-d_1||\alpha_2|}{I_2^2} \right.\right. 
  \left.\left.
  + \frac{|d_1-d_2||\alpha_3|}{I_3^2}\right)\right.\\
  \left.+\beta^2\left(\frac{|d_2-d_3|^2}{I_1^2}
  + \frac{|d_3-d_1|^2}{I_2^2}
  + \frac{|d_1-d_2|^2}{I_3^2}\right)\right]
\end{multline*}
then the previous assumptions imply that~: 
\begin{equation*}
  (e_1^*)^2+(e_2^*)^2+(e_3^*)^2<f(\beta)
\end{equation*}

Therefore, a sufficient condition for the scheme to be bounded is
$f(\beta)\leq\beta$. 
We know that~:
\begin{equation*}
  \frac{|d_2-d_3|}{I_1}=\frac{|d_2-d_3|}{d_2+d_3}\leq1
\end{equation*}
as the $d_i$ are positive.
Then~:
\begin{equation*}
  f(\beta)\leq
  \frac{1}{4(1-\beta)}\left(\Delta t^2\left[\frac{|\alpha_1|^2}{I_1^2}
  + \frac{|\alpha_2|^2}{I_2^2}
  + \frac{|\alpha_3|^2}{I_3^2}\right)\right. 
  \left. + 2\beta \Delta t\left(\frac{|\alpha_1|}{I_1}
  + \frac{|\alpha_2|}{I_2}
  + \frac{|\alpha_3|}{I_3}\right)
  +3\beta^2 \right]
\end{equation*}
Hence, a sufficient condition for $f(\beta)\leq\beta$ to hold is~:
\begin{equation}
  \Delta t^2\left(\frac{|\alpha_1|^2}{I_1^2}+\frac{|\alpha_2|^2}{I_2^2}+\frac{|\alpha_3|^2}{I_3^2}\right) 
 +2\beta
 \Delta t\left(\frac{|\alpha_1|}{I_1}+\frac{|\alpha_2|}{I_2}+\frac{|\alpha_3|}{I_3}\right)
 +7\beta^2-4\beta<0 \label{eqn:condition}
\end{equation}
Let us define~:
\begin{align*}
  B &=
  \frac{|\alpha_1|}{I_1}+\frac{|\alpha_2|}{I_2}+\frac{|\alpha_3|}{I_3}\\
  C &= \frac{|\alpha_1|^2}{I_1^2}+\frac{|\alpha_2|^2}{I_2^2}+\frac{|\alpha_3|^2}{I_3^2}
\end{align*}
A sufficient condition to obtain (\ref{eqn:condition}) is to have
$\Delta t\leq\tilde{\Delta t}$ with~:
\begin{equation*}
  \tilde{\Delta t} = \frac{-2\beta B+\sqrt{4\beta^2B^2-4(7\beta^2-4\beta)C}}{2C}
\end{equation*}
As we supposed that $0<\beta<\frac{1}{2}<\frac{4}{7}$,
$7\beta^2-4\beta<0$. We also know that $B^2\leq3C$ and $C\leq B^2$,
and it follows that~:
\begin{equation*}
  \tilde{h}\geq\frac{2\sqrt{\frac{\beta-\beta^2}{3}}-\beta}{B}
\end{equation*}

In the end, we have the following lemma~:
\begin{lem}
  Let us choose $0<\beta<\frac{1}{2}$ and $\Delta t>0$ such that~:
  \begin{equation}
    \Delta t\left(\frac{|\alpha_1|}{I_1}+\frac{|\alpha_2|}{I_2}+\frac{|\alpha_3|}{I_3}\right)\leq2\sqrt{\frac{\beta-\beta^2}{3}}-\beta \label{eqn:CFL}
  \end{equation}
  If $(e_1,e_2,e_3)\in\mathcal{B}(0,\sqrt{\beta})$, then
  $(e_1^*,e_2^*,e_3^*)\in\mathcal{B}(\vect{0},\sqrt{\beta})$. 
\end{lem}

\subsection{The iterative scheme is a contraction}

Following the previous subsection, suppose that $(e_1,e_2,e_3)$ and
$(f_1,f_2,f_3)$ are in $\mathcal{B}(0,\sqrt{\beta})$, and let $e_0
= \sqrt{1-e_1^2-e_2^2-e_3^2}$ and $f_0 =
\sqrt{1-f_1^2-f_2^2-f_3^2}$. We define $e^*$ and $f^*$ as before. We
show here that $\lVert e^*-f^*\rVert\leq\rho\lVert e-f\rVert$, with
$0<\rho<1$. 

We compute~:
\begin{equation*}
  e_1^*-f_1^* =
  \frac{(d_2-d_3)}{I_1e_0}\left[(f_2-e_2)\left(\frac{f_3+e_3}{2}\right)
  \right. 
  \left.
  +(f_3-e_3)\left(\frac{f_2+e_2}{2}\right)\right]+\frac{f_0-e_0}{e_0}f_1^* 
\end{equation*}

We then use the fact that $\frac{|d_2-d_3|}{I_1}<1$. As the same type
of results hold with a circular permutation of indices $x$, $y$ and $z$, we let
$\lVert\cdot\rVert$ the euclidian norm in $\mathbb{R}^3$ on
$(e_1,e_2,e_3)$, and we find~:
\begin{multline*}
  \lVert e^*-f^*\rVert^2 \leq 
  2\frac{\left(\frac{f_2+e_2}{2}\right)^2+\left(\frac{f_3+e_3}{2}\right)^2}{e_0^2}(f_1-e_1)^2 
  +
  2\frac{\left(\frac{f_1+e_1}{2}\right)^2+\left(\frac{f_3+e_3}{2}\right)^2}{e_0^2}(f_2-e_2)^2 \\
   +
  2\frac{\left(\frac{f_1+e_1}{2}\right)^2+\left(\frac{f_2+e_2}{2}\right)^2}{e_0^2}(f_3-e_3)^2 
  +
  \frac{4}{e_0^2}(f_2-e_2)(f_3-e_3)\left(\frac{f_2+e_2}{2}\right)\left(\frac{f_3+e_3}{2}\right)\\
   +
  \frac{4}{e_0^2}(f_1-e_1)(f_3-e_3)\left(\frac{f_1+e_1}{2}\right)\left(\frac{f_3+e_3}{2}\right) 
  +
  \frac{4}{e_0^2}(f_1-e_1)(f_2-e_2)\left(\frac{f_1+e_1}{2}\right)\left(\frac{f_2+e_2}{2}\right)
  \\
   +2\frac{(f_1^*)^2+(f_2^*)^2+(f_3^*)^2}{e_0^2}(f_0-e_0)^2
\end{multline*}

Since~:
\begin{equation*}
  \frac{4}{e_0^2}(f_2-e_2)(f_3-e_3)\left(\frac{f_2+e_2}{2}\right)\left(\frac{f_3+e_3}{2}\right) 
  \leq \frac{2}{e_0^2}\left[(f_2-e_2)^2\left(\frac{f_2+e_2}{2}\right)^2 \right.
  \left.+ (f_3-e_3)^2\left(\frac{f_3+e_3}{2}\right)^2\right]
\end{equation*}
we have~:
\begin{equation*}
  \lVert e^*-f^*\rVert^2
  \leq \frac{2}{e_0^2}\left(\lVert\frac{e+f}{2}\rVert^2 \lVert e-f\rVert^2 +
  \lVert f^*\rVert^2 (f_0-e_0)^2 \right) 
\end{equation*}

We also have~:
\begin{align*}
  (f_0-e_0)^2 \leq \frac{\lVert \frac{e+f}{2}\rVert^2}{(\frac{e_0+f_0}{2})^2}\lVert e-f\rVert^2
\end{align*}

In the end, we obtain the upper bound~:
\begin{equation*}
  \lVert e^*-f^*\rVert^2 \leq 2\frac{\lVert
  \frac{e+f}{2}\rVert^2}{e_0^2}\left(1+\frac{\lVert
  f^*\rVert^2}{\left(\frac{e_0+f_0}{2}\right)^2}\right) \lVert e-f\rVert^2
\end{equation*}

If we take the same hypotheses as in the first subsection, that is,
$(e_1,e_2,e_3)\in\mathcal{B}(\vect{0},\sqrt{\beta})$ and
$(f_1,f_2,f_3)\in\mathcal{B}(\vect{0},\sqrt{\beta})$, and $h$ such that
$(e_1^*,e_2^*,e_3^*)\in\mathcal{B}(\vect{0},\sqrt{\beta})$ and
$(f_1^*,f_2^*,f_3^*)\in\mathcal{B}(\vect{0},\sqrt{\beta})$, then due to the
convexity of $\mathcal{B}(\vect{0},\sqrt{\beta})$, we have~:
\begin{equation*}
  \lVert \frac{e+f}{2}\rVert^2<\beta
\end{equation*}
and moreover, as $e_0^2>1-\beta$ et $f_0^2>1-\beta$, then
$\left(\frac{e_0+f_0}{2}\right)^2>1-\beta$.

Then~:
\begin{equation*}
  2\frac{\lVert\frac{e+f}{2}\rVert^2}{e_0^2}\left(1+\frac{\lVert
  f^*\rVert^2}{\left(\frac{e_0+f_0}{2}\right)^2}\right) 
  \leq 2\frac{\beta}{1-\beta}\left(1+\frac{\beta}{1-\beta}\right) = \frac{2\beta}{(1-\beta)^2}
\end{equation*}

In order to have a scheme which is a contraction, it is sufficient to
impose~:
\begin{equation*}
  \frac{2\beta}{(1-\beta)^2} \leq 1
\end{equation*}
As $0<\beta<\frac{1}{2}$, it is sufficient to choose~:
\begin{equation*}
  \beta\leq 2-\sqrt{3}
\end{equation*}

\subsection{Optimization on constant $\beta$}

Optimizing the stability condition (\ref{eqn:CFL}) on $\Delta t$, we obtain
the following optimal value of $\beta$~:
\begin{equation*}
  \beta_{\text{max}} = \frac{7-\sqrt{21}}{14} \approx 0.17
\end{equation*}

\subsection{Conclusion}

If we take the time-step $\Delta t$ such that~:
\begin{equation*}
  \Delta t\left(\frac{|\alpha_1|}{I_1}+\frac{|\alpha_2|}{I_2}+\frac{|\alpha_3|}{I_3}\right)
  \leq 2\sqrt{\frac{\beta_{max}-\beta_{max}^2}{3}}-\beta_{max} \approx 0.26
\end{equation*}
then the iterative scheme starting with $(1,0,0,0)$ converges to the
unique solution of the nonlinear problem in
$\mathcal{B}(\vect{0},\sqrt{\frac{7-\sqrt{21}}{14}})$, and the convergence
speed is geometric with a rate $\rho<1$. In addition,
$\rho<28-6\sqrt{21}\approx0.5$. We thus have 
proved existence and uniqueness of the solution in
$\mathcal{B}(\vect{0},\frac{\sqrt{2}}{2})$. 


\end{document}